\documentclass[review]{elsarticle}
\usepackage{graphicx}
\usepackage[colorlinks]{hyperref}
\hypersetup{pdfauthor={Name}}
\usepackage{lineno}
\usepackage{float}
\usepackage{multicol,multirow}
\usepackage{multienum}
\usepackage{eucal}
\usepackage[T1]{fontenc}
\usepackage[utf8]{inputenc}
\usepackage[english]{babel}
\usepackage{mathtools}
\usepackage{amsmath,amssymb}

\usepackage{graphics}

\usepackage{caption}
\usepackage{subcaption}

\usepackage{xcolor}
\usepackage[hmargin=1in,vmargin=1in]{geometry}

\allowdisplaybreaks

\makeatletter
\newcommand{\subjclass}[2][2020]{%
	\let\@oldtitle\@title%
	\gdef\@title{\@oldtitle\footnotetext{#1 \emph{Mathematics subject classification}: #2}}%
}
\newcommand{\keywords}[1]{%
	\let\@@oldtitle\@title%
	\gdef\@title{\@@oldtitle\footnotetext{\emph{Keywords}: #1}}%
}
\makeatother

\newtheorem{theorem}{Theorem}

\newdefinition{remark}{Remark}
\newproof{proof}{Proof}

\modulolinenumbers[5]

\journal{Appropriate journal to be identified}

\begin{document}
	
	\begin{frontmatter}
		\title{Analysis of Heterogeneous Vehicular Traffic: Using Proportional Densities}
		\tnotetext[t1]{This article is a result of the PhD program funded by Makerere-Sida Bilateral Program under Project 316 "Capacity Building in Mathematics and Its Applications".}
		
		\author[mymainaddress,mysecondaryaddress]{Nanyondo Josephine\corref{mycorrespondingauthor}}
		\cortext[mycorrespondingauthor]{Corresponding author}\ead{jnanyondo@sci.busitema.ac.ug}
		
		\author[mysecondaryaddress]{Henry Kasumba}
		\ead{kasumba@cns.mak.ac.ug}
				
			\address[mymainaddress]{Department of Mathematics, Faculty of Science and Education, Busitema University, P.O Box 236, Tororo, Uganda}			
		\address[mysecondaryaddress]{Department of Mathematics, School of Physical Sciences, Makerere University, P.O Box 7062, Kampala, Uganda}
	
		\begin{abstract}
				An extended multi-class Aw-Rascle (AR) model with pressure term described as a function of area occupancy defined in form of proportional densities is presented. Two vehicle classes that is; cars and motorcycles are considered based on an assumption that proportions of these form total traffic density. Qualitative properties of the proposed equilibrium velocity is established. Conditions under which the proposed model is stable are determine by linear stability analysis. To compute numerical flux, the model is discretized by the original Roe decomposition scheme, where Roe matrix, averaged data variables and wave strengths are explicitly derived. The Roe matrix is shown to be hyperbolic, consistent and conservative. From the numerical results, the effect of motorcycles proportion on the flow of vehicle classes is determined. Results obtained remain within limits therefore, the proposed model is realistic.
		\end{abstract}
	\end{frontmatter}
	\linenumbers
	\section{Introduction}\label{sec:Intro}
	In developing countries, various vehicle classes use same roads at the same time. This is why traffic flow models that combine vehicle classes are necessary in the study of heterogeneous vehicular traffic flow. The traffic flow models are classified as macroscopic, microscopic or mesoscopic. Macroscopic models consider traffic flow as a compressible fluid formed by vehicles \cite{li2012microscopic} and capture general relationship between flow, density and velocity of vehicles \cite{khan2019macroscopic,moller2014introduction}. Microscopic models describe interactions between vehicles within a traffic stream \cite{maerivoet2005transportation, knoop2017introduction}. The microscopic models consider characteristics such as; vehicle lengths, speeds, accelerations, time, space headways, vehicle and engine capabilities and human characteristics that describe the driving behaviors. Mesoscopic models are described by a combination of microscopic and macroscopic characteristics through probability distributions \cite{khan2019macroscopic} but at a high level of aggregation usually by speed-density relations and queuing theory approaches \cite{moller2014introduction}. 
	
	Traffic flow is either homogeneous, heterogeneous, equilibrium or non-equilibrium. In homogeneous traffic, all vehicles follow lane discipline whereas in heterogeneous traffic, vehicles violate lane discipline. Equilibrium flow is characterized by constant velocity and spatial homogeneity. Whereas, non-equilibrium traffic flow is characterized by changes in both velocity and spatial homogeneity \cite{khan2019macroscopic}.
	
	The macroscopic models are adopted because of their low computation complexity. The first macroscopic model was proposed by \cite{lighthill1955kinematic} and \cite{richards1956shock}, given by
	\begin{equation*}
		\dfrac{\partial \rho}{\partial t} + \dfrac{\partial (\rho v)}{\partial x}= 0.
	\end{equation*}
	
	To address the drawbacks of the first order methods in characterizing overtaking and creeping some researchers (see for example, \cite{FanWork2015HMTraffic} and the references therein) used various forms of definitions of the velocity function (in static form) together with the equation for the conservation of mass. However, the exact formulation of the velocity function for each class of vehicles is a challenge. Using the porous flow approach, Nair et al \cite{nair2011porous} formulated a multi-class heterogeneous traffic flow model using first order methods. The same theory is applied also in \cite{Gashaw2018MixedFlowTraffic,Hawkins2018HTrafficPorus} to formulate a model for heterogeneous traffic flow that takes overtaking and creeping into account. Since these models are based on the first order method they may fail to characterize some important traffic phenomena that occurs in heterogeneous traffic flow, like, stop and go wave, hysteresis \cite{knoop2017introduction}, platoon dispersion \cite{mohan2013heterogeneous}.
	
	Second and higher order models combine the dynamics of the velocity function together with the conservation equation of vehicles. For instance Payne \cite{Payne1971,Payne1979FREFLO} has shown that the microscopic car-following model can only be generalized by a second and higher order continuum models.
	Gupta and Dhiman \cite{gupta2014analyses} have shown that a heterogeneous continuum second order model generalizes a kinematic car-following model for a non lane-based system of traffic flow that takes lateral separation into account. However, in their model, Gupta and Dhiman used the same equilibrium velocity for all vehicle classes. In order to characterize the phenomena of creeping (or gap filling), Mohan and Ramadurai \cite{mohan2013heterogeneous,mohan2017HTModel} used lateral area occupancy method to analyze the flow of heterogeneous traffic. However, they ignored the concept of proportional densities and their effect on the flow of vehicle classes. 
	
	In this paper, the effect of motorcycles proportion on the flow of vehicle classes is determined. The subsequent subsections of the paper are organized as follows. Derivation of area occupancy defined in terms of proportional densities and conservative form of the proposed model are shown in Section $2$. In Section $3$, qualitative properties of equilibrium velocity are discussed. A discretization scheme, Quasi-linear form, original Roe numerical scheme in which Roe matrix, averaged data variables, wave strengths and entropy condition are shown in Section $4.$ Stability analysis of the proposed model equations is done in Section 5. Numerical results are shown in Section $6$. Finally, Section $7$ draws conclusion on the paper.
	
\section{Model formulation}\label{sec:ModelFormulation}
	In this section, the Aw-Rascle model 
	is modified by introducing an area occupancy that is expressed in terms of vehicle class proportions and other lateral parameters in the pressure term.
	\subsection{Applying area occupancy for heterogeneous traffic}
	We then introduce area occupancy ($AO$) that is defined in terms of proportional densities and derive it as follows.
	Suppose that the total density, $\rho$ of vehicles that are occupying the road, at any time is given by
	\begin{equation}
		\rho = \rho_m + \rho_c, \label{TotalDensity}
	\end{equation}
	where $\rho_c$ refers to density of cars and $\rho_m$  density of motorcycles.
	Let
	\begin{equation}
		\rho_m = \delta(\rho_m + \rho_c) = \delta\rho,
		\label{Rhom}
	\end{equation}
	where $\delta$ represents the proportion of motorcycles that contribute to the total vehicular density on a given road section.
	Then
	\begin{equation}
		\rho_m = \dfrac{\delta  \rho_c}{1- \delta}.
		\label{Rhom2}
	\end{equation}
	Assuming that the width of  a motorcycle is one third the width of a car and considering $l, w$ to be length and width of vehicles, 
	the equation for area occupancy becomes
	\begin{equation}
		AO = \dfrac{\sum_{i} \rho_i a_i}{W} = \dfrac{\rho_c a_c + \rho_m a_m}{W} = \dfrac{(1 - \delta)\rho a_c + \delta \rho a_m}{W},~ i=c,m
		\label{AOc1}.
	\end{equation}
	Substitution of (\ref{Rhom2}) into (\ref{AOc1}) leads to
	\begin{equation}
		AO = \dfrac{1}{W} \left(\rho_c a_c  + \dfrac{\delta \rho_c a_m }{1- \delta}\right) = \dfrac{\rho_c \left[(1- \delta) a_c + \delta a_m\right]}{W(1 -\delta)},
		\label{AO1}\end{equation}
	where $a_i = l_i \times w_i.$
	Further simplification yields
	\begin{eqnarray}
		AO &=& \dfrac{(1- \delta)\rho_c l_c w_c + \delta \rho_c l_m w_m}{W(1 -\delta)}, ~w_m=\dfrac{1}{3}w_c, \nonumber\\
		&=& \dfrac{\rho_c w_c}{W(1- \delta)}\left((1- \delta)l_c + \dfrac{1}{3} \delta l_m\right), \label{AOc2}\\
		& =& \rho_c \psi_c, \nonumber 
	\end{eqnarray}
	where
	\begin{equation}
		\psi_c = \dfrac{w_c}{W(1- \delta)}\left((1- \delta)l_c + \dfrac{1}{3} \delta l_m\right)
		\label{psic2},
	\end{equation}
	and $AO$ refers to area occupancy expressed in terms of proportional densities. Equivalently, $AO$ can be rewritten as
	\begin{eqnarray}
		AO &=& \dfrac{\rho_m w_c(l_c(1- \delta) + \dfrac{1}{3}\delta l_m)}{W \delta}, \label{AOm2}\\
		&=& \rho_m \psi_m, \nonumber
	\end{eqnarray}
	where
	\begin{equation}
		\psi_m = \dfrac{w_c(l_c(1- \delta) + \dfrac{1}{3}\delta l_m)}{W \delta}.
		\label{psim2}
	\end{equation}
	It can be inferred from (\ref{Rhom2}) that
	\begin{eqnarray}
		\rho_c = \dfrac{\rho_m(1- \delta)}{\delta} = (1 - \delta)\rho.
		\label{Rhoc}
	\end{eqnarray} Consequently the total density in (\ref{TotalDensity}) then becomes
	\begin{equation*}
		\rho = \rho_{m} + \rho_{c}
		= \delta\rho + (1-\delta)\rho.
	\end{equation*}
	Incorporating the two vehicle classes in AR model, we obtain the following proposed model
	\begin{eqnarray}
		\begin{cases}
			\dfrac{\partial \rho_m}{\partial t } + \dfrac{ \partial \left(\rho_m v_m\right)}{\partial x } = 0, \\[1.5ex]
			\dfrac{\partial \left(v_m + p_m\right)}{\partial t} + v_m \dfrac{\partial \left(v_m + p_m\right)}{\partial x} = \dfrac{1}{\tau_m}\left(v_{em} - v_m\right), \\[1.5ex]
			\dfrac{\partial \rho_c}{\partial t } + \dfrac{ \partial \left(\rho_c v_c\right)}{\partial x } = 0, \\[1.5ex]
			\dfrac{\partial \left(v_c + p_c\right)}{\partial t} + v_c \dfrac{\partial \left(v_c + p_c\right)}{\partial x} = \dfrac{1}{\tau_c}\left(v_{ec} - v_c\right),
		\end{cases}
		\label{researchmodel1} \end{eqnarray}
	where the respective pressure terms are given by
	\begin{equation*}
		p_m = \left(\psi_m \rho_m\right)^{\gamma_m} ~\text{and} ~p_c = \left(\psi_c \rho_c\right)^{\gamma_c}.
	\end{equation*} Here $\psi_m,~ \psi_c$ are parameters defined as in (\ref{psic2}), (\ref{psim2}) and $\gamma_m,~ \gamma_c$ represent the relaxation time of motorcycles and cars, respectively. 
	
	We can express (\ref{researchmodel1}) in conservative form given by 
	\begin{equation}
		\begin{cases}
			\dfrac{\partial \rho_m}{\partial t} + \dfrac{\partial \left(\rho_m v_m\right)}{\partial x} = 0, \\[1.5ex]
			\dfrac{\partial \left(\rho_m \left(v_m + p_m\right)\right)}{\partial t} + \dfrac{\partial \left(\rho_m v_m \left(v_m + p_m\right)\right)}{\partial x} = \dfrac{\rho_m}{\tau_m}\left(v_{em} - v_m\right) , \\[1.5ex]
			\dfrac{\partial \rho_c}{\partial t} + \dfrac{\partial \left(\rho_c v_c\right)}{\partial x} = 0, \\[1.5ex]
			\dfrac{\partial \left(\rho_c \left(v_c + p_c\right)\right)}{\partial t} + \dfrac{\partial \left(\rho_c v_c \left(v_c + p_c\right)\right)}{\partial x} = \dfrac{\rho_c}{\tau_c}\left(v_{ec}  - v_c\right).
		\end{cases} \label{conserving4} \end{equation} 
	where the conserved variables are $\rho_c, ~ \rho_c\left(v_c + p_c\right), ~ \rho_m$ and $\rho_m\left(v_m + p_m\right),$ \\
	Considering terms in system (\ref{conserving4}) with subscript $c,$ let the generalized momentum be $X_c,$ such that
	\begin{equation}
		X_c = \rho_c\left(v_c + p_c\right).
		\label{Xc1}\end{equation}
	Rearranging terms in (\ref{Xc1}) we obtain
	\begin{eqnarray*}
		\dfrac{X_c}{\rho_c} = \left(v_c + p_c\right),~\text{and}~
		\rho_c v_c = X_c - \rho_c p_c.
	\end{eqnarray*}
	This implies
	\begin{equation*}
		\rho_c v_c \left(v_c + p_c\right) = \left(X_c - \rho_c p_c\right)\dfrac{X_c}{\rho_c},
	\end{equation*} which simplifies into
	\begin{equation*}\dfrac{X_c^2}{\rho_c} - p_cX_c.
	\end{equation*}
	The same is done on terms with subscript $m$ to obtain the following conserved model equations.
	\begin{equation}
		\begin{cases}
			\dfrac{\partial \rho_m}{\partial t} + \dfrac{ \partial \left(X_m - \rho_m p_m\right)}{\partial x} = 0, \\[1.5ex]
			\dfrac{\partial X_m}{\partial t} + \dfrac{\partial \left(\dfrac{X_m^2}{\rho_m} - p_mX_m\right)}{\partial x} = \dfrac{\rho_m}{\tau_m}\left(v_{em} - v_m\right), \\[1.5ex]
			\dfrac{\partial \rho_c}{\partial t} + \dfrac{ \partial \left(X_c - \rho_c p_c\right)}{\partial x} = 0, \\[1.5ex]
			\dfrac{\partial X_c}{\partial t} + \dfrac{\partial \left(\dfrac{X_c^2}{\rho_c} - p_cX_c\right)}{\partial x} = \dfrac{\rho_c}{\tau_c}\left(v_{ec}  - v_c\right).
		\end{cases} \label{conserving5} \end{equation}
	In vector form, system (\ref{conserving5}) becomes
	\begin{equation}
		\dfrac{\partial U}{\partial t} + \dfrac{\partial f(U)}{\partial x} = S(U),
		\label{VectorForm1} \end{equation}
	where
	\begin{equation*}
		U = \begin{bmatrix}
			\rho_m \\[1.5ex]
			X_m \\[1.5ex]
			\rho_c \\[1.5ex]
			X_c\end{bmatrix}, ~
		f(U) =
		\begin{bmatrix}
			X_m - \rho_m p_m \\[1.5ex]
			\dfrac{X_m^2}{\rho_m} - p_mX_m\\[1.5ex]
			X_c - \rho_c p_c \\[1.5ex]
			\dfrac{X_c^2}{\rho_c} - p_cX_c
		\end{bmatrix},~ \\
		S(U) =
		\begin{bmatrix}
			0 \\[1.5ex]
			\dfrac{\rho_m}{\tau_m}\left(v_{em} - v_m\right)\\[1.5ex]
			0 \\[1.5ex]
			\dfrac{\rho_c}{\tau_c}\left(v_{ec}  - v_c\right)
		\end{bmatrix}.
		\label{VectorForm2}
	\end{equation*}
\section{Qualitative properties of Equilibrium velocity}\label{sec:EquilibriumVelocity}
We adopt a modified version of the Greenshield's equilibrium velocity 
\begin{eqnarray*}
	\begin{cases}
		v^e_m = v_{m \max}\left(1 - \dfrac{AO_m}{AO_m^{\max}}\right),~ \text{if}~ AO_m \leq AO_m^{\max},\\[1.5ex] 
		v^e_c = v_{c \max}\left(1 - \dfrac{AO_c}{AO_c^{\max}}\right),~ \text{if}~ AO_c \leq AO_c^{\max}, \\[1.5ex]
		v^e_i = 0 ~ \text{otherwise,} ~i=m,c,
	\end{cases}
\end{eqnarray*}  where $v_{m\max},~ v_{c\max},~ AO_m^{\max},~ AO_c^{\max}$ are the respective maximum velocity and maximum area occupancy of the vehicle classes. Velocities were defined in terms of $AO$ since it is a standard measure of heterogeneous traffic concentration that lacks in lane discipline. 
Note that by writing $AO$ in the form $AO_i =  \psi_i \rho_i,~ i = m,~ c,$ and using definitions 
\begin{equation*}
	w_m = \dfrac{w_c}{3} ~\text{and}~ a_i = l_i*w_i,~ i =m, c,
\end{equation*}
in (\ref{AOm2}) gives 
\begin{equation}
	v^e_m = v_{m \max}\left(1 - \dfrac{\rho_m\left[a_c\left(1 - \delta\right) + \delta a_m\right]}{0.85W \delta}\right). \label{vmeExpressedInFormRhom2}
\end{equation} 

Similarly, using (\ref{AOc2}) yields
\begin{equation}
	v^e_c = v_{c \max}\left(1 - \dfrac{\rho_c\left[a_c\left(1 - \delta\right) + \delta a_m\right]}{ 0.74W (1 - \delta)}\right). \label{vceExpressedFormRhoc2}
\end{equation}
The dependence of equilibrium speed of the flow, $v^e_i,$ on the area occupied by all vehicles is shown in (\ref{dvem2}) and (\ref{dvec2}) by
taking derivatives of $d v^e_c,$ with respect to $AO$ to obtain:
\begin{eqnarray}
	\dfrac{d v^e_m}{d AO_m} &=& - \dfrac{v_{m \max}}{AO_m^{\max}}, \label{dvem2}\\
	\dfrac{d v^e_c}{d AO_c} &=& - \dfrac{v_{c \max}}{AO_c^{\max}} . \label{dvec2}
\end{eqnarray}
From the derivatives above observe that the speed of traffic flow decreases as the area occupied by vehicles increases.

\begin{figure}[H]
	\centering
	\begin{subfigure}[t]{0.45\linewidth}
		\centering 
		\includegraphics[height=1.6in]{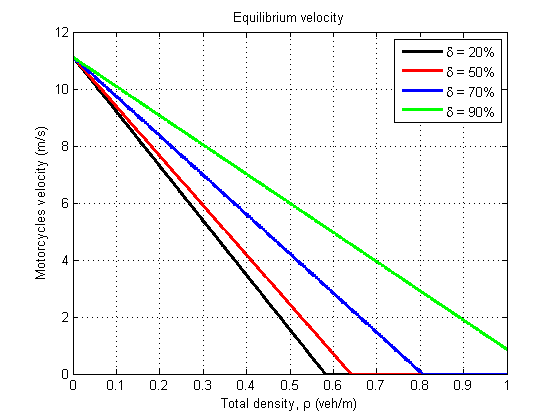}
		\caption{Equilibrium velocity of motorcycles vs total density}\label{EqulSpeedMotorcycle}
	\end{subfigure}
	\begin{subfigure}[t]{0.45\linewidth}
		\centering
		\includegraphics[height=1.6in]{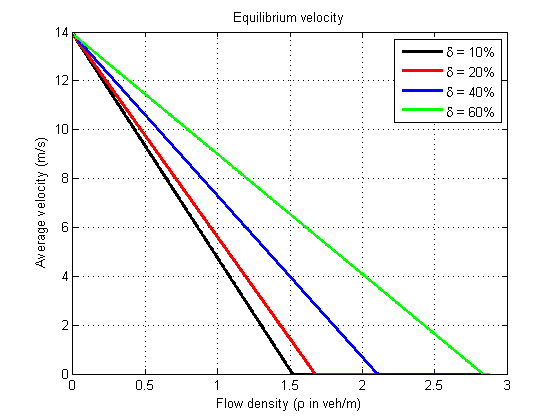}
		\caption{Equilibrium velocity of cars vs total density}\label{EqulSpeedCar}
	\end{subfigure}
	\begin{subfigure}[t]{0.45\linewidth}
		\centering
		\includegraphics[height=1.6in]{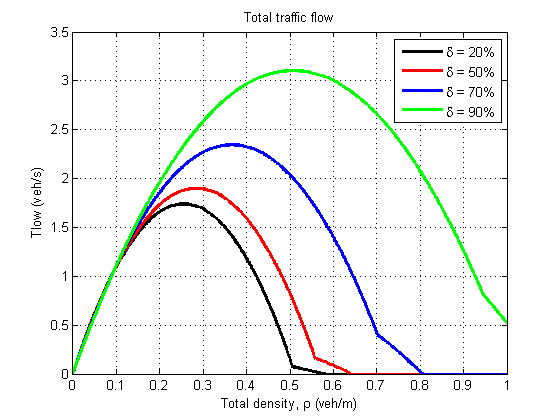}
	\noindent \caption{Total flow vs total density}\label{TotalTrafficFlow}
	\end{subfigure}
	\noindent \caption{Equilibrium velocities and total flow of the two vehicle classes versus total density, with varying values of motorcycles' proportion, $\delta.$}
	\label{EqmVelocityDensityrelation}
	\vspace{-0.1cm}
\end{figure}

By taking derivatives with respect to the motorcycles proportion $\delta,$ we show how the speed of each vehicle class is affected by the variation in the proportion. We make use of \ref{vmeExpressedInFormRhom2} to derive
\begin{eqnarray}
	\dfrac{d v^e_m}{d \delta} 
	&=& \dfrac{\rho_m a_c v_{m\max}}{0.85 \delta^2 W}, \label{dvemdRho2}
\end{eqnarray}
It is concluded from (\ref{dvemdRho2}) that motorcycles' equilibrium speed increases with of their proportion.
For the case of cars, \ref{vceExpressedFormRhoc2}  is used in deriving 
\begin{eqnarray}
	\dfrac{d v^e_c}{d \delta}
	&=& - \dfrac{a_m \rho_{c}v_{c\max}}{W(1- \delta)^2}. \label{dvecdRho2}
\end{eqnarray}
It can be deduced from (\ref{dvemdRho2}) that an increase in motorcycles proportion leads to a decrease in velocity of cars and vise versa.
Next, we look at the dynamics of equilibrium velocity when density is varied. That is achieved through derivatives  
\begin{equation}
	\dfrac{d v_m^e}{d \rho} =  -v_{m \max}\left(\dfrac{\left[a_c\left(1 - \delta\right) + \delta a_m\right]}{0.85W}\right), \label{dvmedrho2}	 
\end{equation}
and
\begin{equation}
	\dfrac{d v_c^e}{d \rho} = -v_{c \max}\left(\dfrac{\left[a_c\left(1 - \delta\right) + \delta a_m\right]}{0.74W}\right), \label{dvcedrho2}	 
\end{equation} 
where $\rho_{m} = \delta \rho,~ \text{and}~ \rho_{c} = (1 - \delta)\rho.$ It is observed from (\ref{dvmedrho2}) and (\ref{dvcedrho2}) that equilibrium velocity of vehicles decreases with increase of density. This behavior is further evidenced in Figure \ref{EqmVelocityDensityrelation}. In Figures \ref{EqulSpeedMotorcycle} and \ref{EqulSpeedCar}, cars are observed to stop much faster than motorcycles irrespective of motorcycles proportion. When the proportion is very large, motorcycles keep moving even when cars have stopped. It is because they can ably maneuver or filter through spaces unfilled by cars. Figure \ref{TotalTrafficFlow} predicts jam density reached at much earlier if cars dominate the road than when motorcycles do. It can be concluded that total flow of vehicles increases with increase of motorcycles proportion.

\section{Discretization scheme}\label{sec:RoeDecomposition}
The proposed model (\ref{VectorForm1}) in a space-time domain $\displaystyle{[0,J] \times [0,N]}$ is considered. Solving it would mean evolving the solution $U(x,t)$ in time starting from the initial condition $U(x, 0)$ and subject to boundary conditions. The space and time grids are divided into $J$ and $N$ equally spaced grid points $x_j = j \Delta x,~ j = 0, \ldots, J$ and $t^n = n \Delta t,~ n = 0 , \ldots, N,$ of length $\Delta x$ and $\Delta t,$ respectively. Approximate or discrete values of $U(x,t),$ for the data given by (\ref{VectorForm1}) at a time step $n$ and spatial position $j,$ are represented by $U_j^n \equiv U \left(j \Delta x,n \Delta t \right) \equiv U \left(x_j, t^n \right).$ Next, periodic conditions are applied to the left and right boundary points $U_0^n$ and $U_J^n,$ respectively. Since at any time step $n,$ the numerical solutions are obtained as either initial data values or values computed in a previous time step, a conservative finite-difference scheme is then applied to (\ref{VectorForm1}). Hence, the solution $U_j^{n+1},$ at the next time step is obtained as
\begin{eqnarray}
	U_j^{n+1} &=& U_j^n - \dfrac{\Delta t}{\Delta x}\left(F(U_j^n, U_{j+1}^n) - F(U_{j-1}^n, U_j^n)\right) + \Delta t S(U_i^n),\label{DifferenceEquationwithSourceterm}
\end{eqnarray}
where $F(U_j^n, U_{j+1}^n)$ in (\ref{DifferenceEquationwithSourceterm}) denotes Roe's numerical flux \citep{toro2013riemann} and is given by
\begin{eqnarray}
	F(U_j^n, U_{j+1}^n) = \dfrac{1}{2} \left(F(U_j^n) + F(U_{j+1}^n)\right) - \dfrac{1}{2} \left(B(U^n_{j+1/2}) \left(U^n_{j+1} - U^n_j \right)\right),
	\label{FluxChangeAtBoundary3}
\end{eqnarray}
where $B(U^n_{j+1/2})$ represents the averaged Jacobian matrix of $f(U)$. For each time step, scheme (\ref{DifferenceEquationwithSourceterm}) is implemented in two intermediate steps. Firstly, one finds
\begin{equation}
	U_j^{*} = U_j^n - \dfrac{\Delta t}{\Delta x}\left(F(U_j^n, U_{j+1}^n) - F(U_{j-1}^n, U_j^n)\right),
	\label{SolutionwithoutSourceterm}
\end{equation}
obtained by ignoring the non-homogeneous part.
Secondly, the solution $U_j^{n+1}$ at time step $t^{n+1}$ is found by adding the non-homogeneous part to (\ref{SolutionwithoutSourceterm}). That is to say;
\begin{equation*}
	U_j^{n+1} = U_j^{*} + \Delta t S(U_j^n).
	\label{SolutionwithSourceterm}
\end{equation*}
Below we describe the Roe's numerical scheme that is utilized in this work.
\subsection{Roe decomposition}
The Roe decomposition scheme is used because it is capable of capturing abrupt changes or discontinuities in the data variables and providing accurate numerical solutions for traffic flow models \cite{khan2019macroscopic}. Roe's scheme is employed purposely to transform the proposed model (\ref{researchmodel1}) into quasi-linear form and linearize locally by approximating the Jacobian matrix, $B(U)$ with Roe averages and in every time step, the procedure is repeated \citep{kachroo2008pedestrian}. The resulting system can then approximate speed found from the eigenvalues of the averaged Jacobian matrix, as described below \citep{kachroo2008pedestrian}.

\subsubsection{Quasi-linear form of model equations}

The model equations are expressed in Quasi-linear form which is suitable for carrying out linearization. The homogeneous part of system (\ref{conserving5}) can be rewritten in quasi-linear form as
\begin{equation*}
	\dfrac{\partial U }{\partial t} + B(U)\dfrac{\partial U}{\partial x} = 0,
	\label{quasilinearform1}
\end{equation*}
where
\begin{eqnarray}
	B(U) &=& \dfrac{\partial f}{\partial U}, \nonumber \\
	&=&
	\begin{bmatrix}
		-(\gamma_m + 1)p_m & 1 & 0 & 0\\
		B_{21} & \dfrac{2X_m}{\rho_m} - p_m & 0 & 0\\
		0 & 0 & -(\gamma_c + 1)p_c & 1 \\
		0 & 0 & - \left(\dfrac{X_c^2}{\rho_c ^2} + \dfrac{\gamma_c p_cX_c}{\rho_c}\right) & \dfrac{2X_c}{\rho_c} - p_c
	\end{bmatrix} \label{ExactJacobianMatrix}
\end{eqnarray}
is the Jacobian matrix with $B_{21}=-\left(\dfrac{X_m^2}{\rho_m ^2} + \dfrac{\gamma_m p_mX_m}{\rho_m}\right).$ Next, we require eigenvalues and eigenvectors of (\ref{conserving5}) to be applied in the Roe numerical flux and entropy fix. Note that the eigenvalues are obtained from the characteristic equations of $B(U)$ and are given as follows:

\begin{equation}
	\lambda_1(U) = v_m,~
	\lambda_2(U) = v_m - \gamma_m p_m,~
	\lambda_3(U) = v_c,~
	\lambda_4(U) = v_c - \gamma_c p_c.
	\label{Eigenvalues}\end{equation}	
The right eigenvectors corresponding to (\ref{Eigenvalues}) are obtained as:
\begin{eqnarray}
	r_1 = \left(\begin{array}{c}
		1\\ (1+\gamma_m)p_m + v_m\\ 0\\ 0 \end{array}\right),~
	r_2 &=& \left(\begin{array}{c} 1\\ v_m + p_m\\ 0\\ 0
	\end{array}\right),~ \nonumber\\
	r_3 = \left(\begin{array}{c}
		0\\ 0\\ 1\\ (1+\gamma_c)p_c + v_c \end{array}\right),~
	r_4 &=& \left(\begin{array}{c} 0\\ 0\\ 1\\ v_c + p_c
	\end{array}\right).
	\label{Eigenvectors}
\end{eqnarray}
We then proceed by applying the original Roe method to (\ref{conserving5}) in order to compute Roe averages. The detailed procedure of the method is described in the following steps \citep{toro2013riemann}. Firstly, a parameter vector $Z,$ is introduced such that $U$ and $f(U)$ of (\ref{VectorForm1}) can be rewritten as \begin{equation*}
	U = U(Z),~ f(U) = f(U(Z)).
\end{equation*}
Next, the jumps
\begin{equation}
	\Delta U = U_r - U_l ,~ \Delta f = f(U_r) - f(U_l),
	\label{ChangeinUF1}
\end{equation}
are expressed in terms of the change
\begin{equation*}
	\Delta Z = Z_r - Z_l.
	\label{ChangeinZ}
\end{equation*}
Letting the parameter vector
\begin{equation*}
	Z =\left[\begin{array}c
		\alpha_1 \\
		\beta_1 \\
		\alpha_2 \\
		\beta_2 \end{array} \right]
	= \left[\begin{array}c
		\sqrt{\rho_m} \\
		\dfrac{X_m}{\sqrt{\rho_m}}\\
		\sqrt{\rho_c} \\
		\dfrac{X_c}{\sqrt{\rho_c}}
	\end{array} \right],
\end{equation*}
implies that $\alpha_1^2 = \rho_m,~  X_m = \alpha_1 \beta_1,~ \alpha_2^2 = \rho_c ~\text{and}~ X_c = \alpha_2 \beta_2.$ So, the vectors 
$	U = U(Z)~\text{and}~ f = f(Z)$ become
\begin{equation*}
	U = \left[\begin{array}c
		\alpha_1^2 \\
		\alpha_1 \beta_1\\
		\alpha_2^2\\
		\alpha_2 \beta_2
	\end{array}\right],~
	f = \left[\begin{array}c
		\alpha_1 \beta_1 - \psi^\gamma \alpha_1^{2(1+ \gamma_m)} \\
		\beta_1^2 - \psi^\gamma \beta_1 \alpha_1^{(1+2 \gamma_m)}\\
		\alpha_2 \beta_2 - \psi^\gamma \alpha_2^{2(1+ \gamma_c)}\\
		\beta_2^2 - \psi^\gamma \beta_2 \alpha_2^{(1+2 \gamma_c)}
	\end{array}\right].
\end{equation*}
The averages of vector $Z$ then become $\overline{Z} = [\overline{\alpha}_1,~\overline{\beta}_1,~\overline{\alpha}_2,~\overline{\beta}_2] ^T.$ These are found by simple arithmetic averaging and obtained as:
\begin{eqnarray*}
	\overline{Z} = \dfrac{1}{2}\left[\begin{array}c
		Z_l + Z_r \end{array} \right]
	=
	\dfrac{1}{2}\left[\begin{array}c
		\alpha_{1l}+\alpha_{1r} \\
		\beta_{1l}+\beta_{1r}\\
		\alpha_{2l}+\alpha_{2r}\\
		\beta_{2l}+\beta_{2r}
	\end{array} \right].
	\label{AveragesofZ}
\end{eqnarray*}
Following work by \cite{toro2013riemann}, two matrices $\overline{B} = \overline{B}(\overline{Z})$ and $\overline{C} = \overline{C}(\overline{Z})$ are computed such that the jumps
$\Delta U$ and $\Delta f$ in (\ref{ChangeinUF1}) are expressed in terms of the jump $\Delta Z.$ That is;
\begin{equation}
	\Delta U = \overline{B}\Delta Z~ \text{and}~ \Delta f = \overline{C}\Delta Z.
	\label{ChangeinUF2}
\end{equation}
Using (\ref{ChangeinUF2}) leads to
$\Delta f = \overline{C}\overline{B}^{-1}\Delta U,$ and comparing it with property (C) of Theorem \ref{RoeMatrixTheorem}, equation (\ref{PropertyCofTheorem}), results into the Roe
matrix $\overline{A} =\overline{C}\overline{B}^{-1}.$ Matrices $\overline{B}$ and $\overline{C}$ satisfying (\ref{ChangeinUF2}) are:
\begin{eqnarray*}
	\overline{B} = \dfrac{dU}{dZ}_{\mid{\overline{Z}}}
	= \left( \begin{array}{cccc}
		2 \overline{\alpha}_1&0&0&0\\
		\overline{\beta}_1& \overline{\alpha}_1&0&0\\
		0&0&2 \overline{\alpha}_2&0\\
		0&0& \overline{\beta}_2& \overline{\alpha}_2
	\end{array} \right),
\end{eqnarray*}
whose inverse is computed by applying the usual definition\\ $\overline{B}^{-1} = \dfrac{1}{Det(\overline{B})} \times Adj(\overline{B})$ and 
\begin{eqnarray*}
	\overline{C} &=& f'(U(Z))\mid_{\overline{Z}},
	\nonumber \\
	&=& \left(\begin{array}{cccc}
		\overline{\beta}_1-2(\gamma_m+1)\psi_m^{\gamma_m} \overline{\alpha}_1^{(1+2 \gamma_m)} & \overline{\alpha}_1 &0 & 0\\
		-(1+2 \gamma_m)\psi_m^{\gamma_m} \overline{\beta}_1 \overline{\alpha}_1^{2 \gamma_m} & \overline{C}_{22} & 0 & 0\\
		0 & 0& \overline{\beta}_2-2(\gamma_c+1)\psi_c^{\gamma_c} \overline{\alpha}_2^{(1+2 \gamma_c)}& \overline{\alpha}_2\\
		0 & 0& -(1+2 \gamma_c)\psi_c^{\gamma_c} \overline{\beta}_2 \overline{\alpha}_2^{2 \gamma_c} & \overline{C}_{44}
	\end{array}\right), \label{Cbar}
\end{eqnarray*} 
where $\overline{C}_{22}=2 \overline{\beta}_1-\psi^\gamma_m \overline{\alpha}_1^{(1+2 \gamma_m)},~\overline{C}_{44}=2 \overline{\beta}_2-\psi^\gamma_c \overline{\alpha}_2^{(1+2 \gamma_c)}.$ Therefore, the required Roe averaged matrix is obtained as
\begin{equation}
	\overline{A} = \left(\begin{array}{cccc}
		\overline{A}_{11} &1 &0 &0\\
		\overline{A}_{21} &2\left(\dfrac{\overline{\beta}_1}{\overline{\alpha}_1}\right) - \psi_m^{\gamma_m}\overline{\alpha}_1^{2\gamma_m} &0 &0\\
		0& 0& -(1+\gamma_c)\psi_c^{\gamma_c} \overline{\alpha}_2^{2\gamma_c} &1\\
		0& 0& \overline{A}_{43}&2\left(\dfrac{\overline{\beta}_2}{\overline{\alpha}_2}\right) - \psi_c^{\gamma_c}\overline{\alpha}_2^{2\gamma_c}
	\end{array}\right), \label{RoeMatrix}
\end{equation}
where $\overline{A}_{11} = -(1+\gamma_m)\psi_m^{\gamma_m} \overline{\alpha}_1^{2\gamma_m},~
\overline{A}_{21} = -\left(\gamma_m\psi^{\gamma_m}\overline{\alpha}_1^{2\gamma_m}\dfrac{\overline{\beta}_1}{\overline{\alpha}_1}  + \dfrac{\overline{\beta}_1^2}{\overline{\alpha}_1^2}\right),\\
\overline{A}_{43} = -\left(\gamma_c\psi^{\gamma_c}\overline{\alpha}_2^{2\gamma_c}\dfrac{\overline{\beta}_2}{\overline{\alpha}_2}  + \dfrac{\overline{\beta}_2^2}{\overline{\alpha}_2^2}\right)~ \text{and}~  \dfrac{\overline{\beta}_1}{\overline{\alpha}_1},$ $\dfrac{\overline{\beta}_2}{\overline{\alpha}_2}$ represent the sum of Roe averaged velocity and pressure of motorcycles and cars, respectively. These are given by
\begin{eqnarray}
	\dfrac{\overline{\beta}_1}{\overline{\alpha}_1} = \dfrac{\dfrac{X_{mr}}{\sqrt{\rho_{mr}}} + \dfrac{X_{ml}}{\sqrt{\rho_{ml}}}}{\sqrt{\rho_{mr}} + \sqrt{\rho_{ml}}} = \overline{v}_m+\overline{p}_m, ~
	\dfrac{\overline{\beta}_2}{\overline{\alpha}_2} = \dfrac{\dfrac{X_{cr}}{\sqrt{\rho_{cr}}} + \dfrac{X_{cl}}{\sqrt{\rho_{cl}}}}{\sqrt{\rho_{cr}} + \sqrt{\rho_{cl}}} = \overline{v}_c+\overline{p}_c.
	\label{RoeAverages}
\end{eqnarray}

Using the Roe averaged velocity, $\overline{v}$ and pressure, $\overline{p},$ the data variables of the proposed model can be approximated over a given road section.
Comparing (\ref{RoeMatrix}) with the matrix in (\ref{ExactJacobianMatrix}), no averaged $\overline{\rho}$ is required. Having derived the Roe matrix, $\overline{A}$ one ensures that it satisfies properties in Theorem \ref{RoeMatrixTheorem} below.

\begin{theorem}[Toro \cite{toro2013riemann}]\label{RoeMatrixTheorem}
	The proposed model is hyperbolic if the Roe Jacobian matrix $\overline{A}$ satisfies the three properties namely:
	\begin{enumerate}
		\item[(A)] The proposed model is hyperbolic if $\overline{A}$ has real eigenvalues $\overline{\lambda}_i = \overline{\lambda}_i(U_l,~U_r),$ where $\overline{\lambda}_1 \leq \overline{\lambda}_2 \leq \cdots \leq \overline{\lambda}_m$ and $\overline{r}_1,~ \overline{r}_2, \ldots, \overline{r}_m$ linearly independent right eigenvectors.
		\item[(B)] Consistency with the exact Jacobian
		$\overline{A}(U,U) = B(U).$
		\item[(C)] Conservation across discontinuities
		\begin{equation}
			f(U_r) - f(U_l) = \overline{A} (U_r - U_l).
			\label{PropertyCofTheorem}\end{equation}
	\end{enumerate}
\end{theorem}
Recalling the Jacobian matrix, eigenvalues and corresponding right eigenvectors of the proposed model given in (\ref{ExactJacobianMatrix}) - (\ref{Eigenvectors}), the averaged eigenvalues are deduced as
\begin{equation*}
	\overline{\lambda}_1(\overline{U}) = \overline{v}_m- \gamma_m\overline{p}_m,~
	\overline{\lambda}_2(\overline{U}) = \overline{v}_m,~
	\overline{\lambda}_3(\overline{U}) =\overline{v}_c- \gamma_c\overline{p}_c,~
	\overline{\lambda}_4(\overline{U}) = \overline{v}_c.
	\label{AveragedEigenvalues}
\end{equation*}
These are shown to be real. Their corresponding averaged right eigenvectors
\begin{eqnarray*}
	\overline{r}_1 &=& \left(\begin{array}{c} 1\\ \overline{v}_m +\overline{p}_m\\ 0\\ 0
	\end{array}\right),~
	\overline{r}_2  = \left(\begin{array}{c}
		1\\ \overline{v}_m +(1+\gamma_m)\overline{p}_m\\ 0\\ 0 \end{array}\right),~ \nonumber \\
	\overline{r}_3 &=& \left(\begin{array}{c} 0\\ 0\\ 1\\ \overline{v}_c +\overline{p}_c
	\end{array}\right),~
	\overline{r}_4 = \left(\begin{array}{c}
		0\\ 0 \\ 1\\ \overline{v}_c +(1+\gamma_c)\overline{p}_c\end{array}\right),
	\label{AveragedEigenvectors}
\end{eqnarray*}
are checked to be linearly independent. Hence, property (A) of Theorem \ref{RoeMatrixTheorem} is satisfied. Property (B), is proved by evaluating the derived Roe matrix (\ref{RoeMatrix}) at $U = \left[ \rho_m,~ X_m,~ \rho_c,~ X_c \right]^T.$ Property $(C)$ is proved by substituting the conserved variables into (\ref{RoeMatrix}). To find the wave strengths, $s_i$ we project the jump, $\Delta U = U_r - U_l$ onto the right eigenvectors and obtain
\begin{eqnarray*}
	\Delta U &=& U_r - U_l = \sum_{i=1}^{4} s_i\overline{r}_i,\\
	\Delta \rho_m &=& \rho_{mr} - \rho_{ml} = s_1 + s_2,\\
	\Delta X_m &=& X_{mr} - X_{ml} = \left(\overline{v}_m+\overline{p}_m\right)s_1 + \left(\overline{v}_m+\left(1+\gamma_m\right)\overline{p}_m\right)s_2,\\
	\Delta \rho_c &=& \rho_{cr} - \rho_{cl} = s_3 + s_4,\\
	\Delta X_c &=& X_{cr} - X_{cl}	= \left(\overline{v}_c+\overline{p}_c\right)s_3 + \left(\overline{v}_c+\left(1+\gamma_c\right)\overline{p}_c\right)s_4.
\end{eqnarray*}
The solutions $s_i$ are verified to be
\begin{eqnarray}\begin{cases}
		s_1 = \Delta \rho_{m} - \dfrac{\Delta X_m - \left(\overline{v}_m + \overline{p}_m\right)\Delta \rho_m}{\gamma_m\overline{p}_m}, \\
		s_2 = \dfrac{\Delta X_m - \left(\overline{v}_m + \overline{p}_m\right)\Delta \rho_m}{\gamma_m\overline{p}_m}, \\
		s_3 = \Delta \rho_{c} - \dfrac{\Delta X_c - \left(\overline{v}_c + \overline{p}_c\right)\Delta \rho_c}{\gamma_c\overline{p}_c}, \\
		s_4 = \dfrac{\Delta X_c - \left(\overline{v}_c + \overline{p}_c\right)\Delta \rho_c}{\gamma_c\overline{p}_c},
	\end{cases} \label{WaveStrengths}
\end{eqnarray}
with the obvious definitions $\Delta \rho_c = \rho_{cr} - \rho_{cl},~ \Delta \rho_m = \rho_{mr} - \rho_{ml},~ \Delta X_c = X_{cr} - X_{cl},~ \Delta X_m = X_{mr} - X_{ml}$. The expressions of $\Delta X_m,~ \Delta X_c$ are given by
\begin{eqnarray*}
	\Delta X_m &=& \Delta\left(\rho_m\left(v_m + p_m\right)\right)	= \Delta \left(\rho_m v_m\right) + \Delta \left(\rho_m p_m\right),\\
	\Delta X_c &=& \Delta\left(\rho_c\left(v_c + p_c\right)\right)	= \Delta \left(\rho_c v_c\right) + \Delta \left(\rho_c p_c\right),
\end{eqnarray*}
where
\begin{eqnarray*}
	\Delta \left(\rho_m v_m\right) &=& \overline{\rho}_m\Delta v_m + \overline{v}_m\Delta \rho_m + O_1(\Delta^2),\\
	\Delta \left(\rho_m p_m\right) &=& \overline{\rho}_m\Delta p_m + \overline{p}_m\Delta \rho_m + O_2(\Delta^2),\\
	\Delta \left(\rho_c v_c\right) &=& \overline{\rho}_c\Delta v_c + \overline{v}_c\Delta \rho_c + O_3(\Delta^2),\\
	\Delta \left(\rho_c p_c\right) &=& \overline{\rho}_c\Delta p_c + \overline{p}_c\Delta \rho_c + O_4(\Delta^2),\\
	O_1(\Delta^2) &=& \left(\rho_{mr}-\rho_{ml}\right)\left(v_{mr}-\overline{v}_m\right) -  \left(\rho_{ml}-\overline{\rho}_m\right)\left(v_{ml}-\overline{v}_m\right),\\
	O_2(\Delta^2) &=& \left(\rho_{mr}-\rho_{ml}\right)\left(p_{mr}-\overline{p}_m\right)  -  \left(\rho_{ml}-\overline{\rho}_m\right)\left(p_{ml}-\overline{p}_m\right), \\
	O_3(\Delta^2) &=& \left(\rho_{cr}-\rho_{cl}\right)\left(v_{cr}-\overline{v}_c\right) -  \left(\rho_{cl}-\overline{\rho}_c\right)\left(v_{cl}-\overline{v}_c\right),\\
	O_4(\Delta^2) &=& \left(\rho_{cr}-\rho_{cl}\right)\left(p_{cr}-\overline{p}_c\right) -  \left(\rho_{cl}-\overline{\rho}_c\right)\left(p_{cl}-\overline{p}_c\right).
\end{eqnarray*}
Finally, using (\ref{RoeAverages}) -- (\ref{WaveStrengths}) into any of the expressions (\ref{FluxChangeAtBoundary3}) or
\begin{equation*}
	\overline{F} = \dfrac{1}{2}\left(F_l + F_r\right) - \dfrac{1}{2}\sum_{i=1}^{4}s_i\mid\overline{\lambda}_i\mid\overline{r}_i,
\end{equation*} leads to the sought Roe numerical flux. 

\subsubsection{Entropy fix} \label{Entropy} 

In this subsection, discontinuities that develop particularly at road section boundaries are dissolved by applying an entropy condition to the Roe decomposition. By doing this, a continuous numerical solution is obtained. An entropy solution $e \mid \Gamma \mid e^{-1},$ is then substituted for the Jacobian matrix (\ref{ExactJacobianMatrix}), where 
$ \mid \Gamma \mid = \left[ \overline{\lambda}_1, ~ \overline{\lambda}_2, ~ \ldots, ~ \overline{\lambda}_K\right]$
refers to a diagonal matrix made up of averaged eigenvalues $\overline{\lambda}_k,~ k=1, 2, \ldots, K$ of (\ref{ExactJacobianMatrix}), $e$ is a matrix made up of averaged eigenvectors and $e^{-1}$ an inverse of $e$ \citep{khan2019macroscopic}. Next, the eigenvalues are modified by employing Harten and Hyman entropy fix scheme in order to accurately characterize the numerical solutions \citep{khan2019macroscopic}. Hence, the modified eigenvalues become
\begin{eqnarray*}
	\overline{\lambda}_k = \begin{cases}
		\overline{\delta}_k,~ \mbox{for}~ \mid \lambda_k \mid \leq \overline{\delta}_k,\\
		\mid \overline{\delta}_k \mid,~ \mbox{for}~ \mid \lambda_k \mid \geq \overline{\delta}_k,
	\end{cases}
\end{eqnarray*} where $\overline{\delta}_k = \max \left(0,~ \lambda_{j+1/2} - \lambda_j,~ \lambda_{j+1} - \lambda_{j+1/2}\right).$ This ensures that the eigenvalues $\overline{\delta}_k$ are non negative and have different values at the road section boundaries. However, for abrupt changes at the road section boundaries, $\overline{\delta}_k$ is set to zero. Using similar procedure as in \cite{khan2019macroscopic}, the transformed Jacobian matrix is given by $e \mid \Gamma \mid e^{-1},$ where
\begin{eqnarray*}
	e &=& \left(\begin{array}{cccc} r_1 & r_2 & r_3 & r_4
	\end{array}\right) \mid_{\overline{\rho}_m,\overline{v}_m\overline{\rho}_c,\overline{v}_cm}, \\
	\mid \Gamma \mid &=& \left(\begin{array}{cccc} \mid \overline{\lambda_1}(U) \mid &0 &0 &0\\
		0& \mid \overline{\lambda_2}(U) \mid &0 &0 \\
		0& 0& \mid \overline{\lambda_3}(U) \mid &0 \\
		0& 0& 0& \mid \overline{\lambda_4}(U) \mid \end{array}\right), \\
	e^{-1} &=&  \left(\begin{array}{cccc}
		\dfrac{-(\overline{v}_m+\overline{p}_m)}{\gamma_m\overline{p}_m} & \dfrac{1}{\gamma_m\overline{p}_m} &0 &0 \\
		\dfrac{\gamma_m\overline{p}_m+(\overline{v}_m+\overline{p}_m)}{\gamma_m\overline{p}_m} & \dfrac{-1}{\gamma_m\overline{p}_m} &0 &0\\
		0& 0& \dfrac{-(\overline{v}_c+\overline{p}_c)}{\gamma_c\overline{p}_c}& \dfrac{1}{\gamma_c\overline{p}_c} \\
		0& 0& \dfrac{\gamma_c\overline{p}_c+(\overline{v}_c+\overline{p}_c)}{\gamma \overline{p}_c}& \dfrac{-1}{\gamma \overline{p}_c}
	\end{array}\right), \\
	\overline{p}_m &=& \psi^\gamma \left(\dfrac{\rho_{mr}+\rho_{mr}}{2}\right)^\gamma,~
	\overline{p}_c = \psi^\gamma \left(\dfrac{\rho_{cl}+\rho_{cr}}{2}\right) ^\gamma,
	\label{ApproxJacobianMatrix2} \end{eqnarray*}
The sought Roe numerical flux is finally derived by evaluating $F(U_j), ~ F(U_{j+1})$ and substituting $ e \mid \Gamma \mid e^{-1}$ for $B(U_{j+1/2})$ into (\ref{FluxChangeAtBoundary3}). Hence, the solutions at subsequent time steps $n+1,$ are obtained using the derived flux (\ref{FluxChangeAtBoundary3}) into (\ref{DifferenceEquationwithSourceterm}). Stability of the discretization schem is established by satisfying the Courant, Friedrich and Lewy (CFL) \cite{leveque1992numerical}.	
	
	\section{Stability analysis of model equations}\label{Stability}
	
	In this Section, behavior of the proposed traffic flow model when subjected to small disturbances is determined { following the procedure outlined in} \cite{gupta2014analyses}. To determine the model's stability, initial densities $\rho_{m0},~ \rho_{c0}$ at time $t_0$ and corresponding equilibrium velocities, $v_{m0}, ~ v_{c0} = v_{em}(\psi \rho_{c0}), ~ v_{ec}(\psi \rho_{c0})$ are considered \cite{khan2019macroscopic}. The initial data is then subjected to changes (or perturbations) denoted $\sigma_{\rho_{m}} (x,t), ~ \sigma_{\rho_{c}} (x,t)$ and $\sigma_{vm} (x,t), ~ \sigma_{vc} (x,t),$ and given by
	\begin{equation*}
		\begin{cases}
			\sigma_{\rho m}(x,t) = \rho_m(x,t) - \rho_{m0},\\
			\sigma_{vm} (x,t) = v_m(x,t) - v_{m0}, \\
			\sigma_{\rho c}(x,t) = \rho_c(x,t) - \rho_{c0},\\
			\sigma_{vc} (x,t) = v_c(x,t) - v_{c0}.
		\end{cases}
		\label{ChangesInDensityandAcceleration}
	\end{equation*}
	Conditions under which changes in density and velocity grows, and eventually causes traffic congestion are determined by assuming sufficiently small changes, smooth functions and periodic functions \cite{gupta2014analyses}. Thus, stability of the proposed model is implied by stability of a linear combination of such functions. Hence, the changes are characterized as
	\begin{equation}
		\begin{cases}
			\sigma_{\rho m}(x,t) = \rho_{m0} e^{ikx+rt},\\
			\sigma_{vm} (x,t) = v_{m0} e^{ikx+rt}, \\
			\sigma_{\rho c}(x,t) = \rho_{c0} e^{ikx+rt},\\
			\sigma_{vc} (x,t) = v_{c0} e^{ikx+rt},
		\end{cases}
		\label{CharacterizedChangesInDensityandAcceleration}
	\end{equation}
	where $i = \sqrt{-1},$ $r$ refers to frequency of oscillations, $k$ number of changes caused by spatial evolution and $kx$ the spatial evolution \cite{khan2019macroscopic}. It can be inferred from (\ref{CharacterizedChangesInDensityandAcceleration}) that traffic flow is a periodic function of $kx$ since\\ $e^{ikx} = \cos(kx) + i \sin(kx).$ At time $t,$ changes in densities and velocities of motorcycles and cars are represented by $\rho_{m0}e^{rt}, ~ \rho_{c0}e^{rt},$ and $v_{m0}e^{rt}, ~ v_{c0}e^{rt},$ respectively and $rt$ is the growth rate of the changes. Therefore, the changes with respect to space and time become 
	\begin{equation*}
		\begin{array}{cc}
			\dfrac{\partial \sigma_{\rho c}(x,t)}{\partial x} = ik \rho_{c0} e^{ikx+rt}, & \dfrac{\partial \sigma_{\rho m}(x,t)}{\partial x} = ik \rho_{m0} e^{ikx+rt}, \\
			\dfrac{\partial \sigma_{\rho c}(x,t)}{\partial t} = r \rho_{c0} e^{ikx+rt}, & \dfrac{\partial \sigma_{\rho m}(x,t)}{\partial t} = r \rho_{m0} e^{ikx+rt},\\
			\dfrac{\partial \sigma_{vc}(x,t)}{\partial x} = ik v_{c0} e^{ikx + rt}, & \dfrac{\partial \sigma_{vm}(x,t)}{\partial x} = ik v_{m0} e^{ikx + rt}, \\
			\dfrac{\partial \sigma_{vc}(x,t)}{\partial t} = rv_{c0} e^{ikx + rt}, & \dfrac{\partial \sigma_{vm}(x,t)}{\partial t} = rv_{m0} e^{ikx + rt}.
		\end{array}
		\label{CharacterizedChangesInDensityandAcceleration2}
	\end{equation*}
	Suppose that a sufficiently small and smooth perturbation $(\sigma_{\rho c}, \sigma_{vc}, \sigma_{\rho m}, \sigma_{vm})$ is given to a steady-state solution $\rho_{c0}, v_{ec}(\rho_{c0}), \rho_{m0}, v_{em}(\rho_{m0}).$ We then examine the effect of the parameter $\psi,$ on the
	propagation of small disturbance. Substituting a small change (\ref{CharacterizedChangesInDensityandAcceleration}) into the proposed model, using Taylor's series expansion at $\left(\rho_{c0}, v_{c0}, \rho_{m0}, v_{m0}\right),$ collecting linear terms where necessary and dropping all non-linear terms \cite{gupta2014analyses} the following linearized system is obtained:
	\begin{eqnarray}
		\begin{cases}
			\dfrac{\partial \sigma_{\rho m}}{\partial t } + v_{m0} \dfrac{ \partial \sigma_{\rho m}}{\partial x } + \rho_{m0} \dfrac{ \partial \sigma_{vm}}{\partial x } = 0,\\
			\dfrac{\partial \sigma_{vm}}{\partial t} + v_{m0} \dfrac{\partial \sigma_{vm}}{\partial x} + \dfrac{\partial p_m(\rho_{m0})}{\partial \rho_m}\dfrac{\partial \sigma_{\rho_m}}{\partial t} + v_{m0} \dfrac{\partial p_m(\rho_{m0})}{\partial \rho_m}\dfrac{\partial \sigma_{\rho m}}{\partial x} - A_1= 0,\\
			\dfrac{\partial \sigma_{\rho c}}{\partial t } + v_{c0} \dfrac{ \partial \sigma_{\rho c}}{\partial x } + \rho_{c0} \dfrac{ \partial \sigma_{vc}}{\partial x } = 0,\\
			\dfrac{\partial \sigma_{vc}}{\partial t} + v_{c0} \dfrac{\partial \sigma_{vc}}{\partial x} + \dfrac{\partial p_c(\rho_{c0})}{\partial \rho_c}\dfrac{\partial \sigma_{\rho_c}}{\partial t} + v_{c0} \dfrac{\partial p_c(\rho_{c0})}{\partial \rho_c}\dfrac{\partial \sigma_{\rho c}}{\partial x} - A_2= 0,
		\end{cases}
		\label{Linearized1}
	\end{eqnarray}
	where $A_1 = \dfrac{\psi_m v_{em}'(\psi_m \rho_{m0}) \sigma_{\rho m} - \sigma_{vm}}{\tau_c},~ A_2 = \dfrac{\psi_c v_{ec}'(\psi_c \rho_{c0}) \sigma_{\rho c} - \sigma_{vc}}{\tau_c}.$
	Following \cite{gupta2014analyses}, we assume that
	\begin{equation*}
		\dfrac{\sigma_{\rho c}(x,t)}{\rho_{c0}} \ll 1, ~ \dfrac{\sigma_{\rho m}(x,t)}{\rho_{m0}} \ll 1,~ \dfrac{\sigma_{vc}(x,t)}{v_{c0}} \ll 1, ~ \dfrac{\sigma_{vm}(x,t)}{v_{m0}} \ll 1.
	\end{equation*}
	Hence, system (\ref{Linearized1}) simplifies to
	\begin{eqnarray*}
		\begin{cases}
			r \sigma_{\rho m} + v_{m0} ik \sigma_{\rho m} + \rho_{m0} ik \sigma_{vm} = 0, \\
			r \sigma_{vm} + v_{m0}ik \sigma_{vm} + r \dfrac{\partial p_m(\rho_{c0})}{\partial \rho_m} \sigma_{\rho m} + v_{m0}ik \dfrac{\partial p_m(\rho_{m0})}{\partial \rho_m} \sigma_{\rho m} - A_2 = 0,\\
			r \sigma_{\rho c} + v_{c0} ik \sigma_{\rho c} + \rho_{c0} ik \sigma_{vc} = 0, \\
			r \sigma_{vc} + v_{c0}ik \sigma_{vc} + r \dfrac{\partial p_c(\rho_{c0})}{\partial \rho_c} \sigma_{\rho c} + v_{c0}ik \dfrac{\partial p_c(\rho_{c0})}{\partial \rho_c} \sigma_{\rho c} - A_1 = 0.
		\end{cases}
	\end{eqnarray*}
	The resulting Jacobian matrix, $J$ is given by
	\begin{equation*}
		\begin{bmatrix}
			r+v_{c0}ik & \rho_{c0} ik &0 &0\\
			\left(r+v_{c0}ik \right)\dfrac{\partial p_c(\rho_{c0})}{\partial \rho_c}-
			\dfrac{\psi_c v_{ec}'(\psi_c \rho_{c0})}{\tau_c} & r+v_{c0}ik+\dfrac{1}{\tau_c} &0 &0 \\
			0& 0& r+v_{m0}ik & \rho_{m0} ik \\
			0& 0& J_{43} & J_{44}
		\end{bmatrix},
	\end{equation*} where
	$
	J_{43} = \left(r+v_{m0}ik \right)\dfrac{\partial p_m(\rho_{c0})}{\partial \rho_m}-
	\dfrac{\psi_m v_{em}'(\psi_m \rho_{m0})}{\tau_m},~ J_{44} = r+v_{m0}ik+\dfrac{1}{\tau_m}.
	$
	For non-trivial solutions, $det(J) = 0.$ Thus, the solutions of $J(\sigma_{\rho c}, \sigma_{vc}, \sigma_{\rho m}, \sigma_{vm}) = 0$ are denoted $r_{1 \pm},~ r_{2 \pm},$ respectively and given by
	\begin{equation*}
		\dfrac{\dfrac{\partial p_m}{\partial \rho_m}ik \rho_{m0} \tau_m -2ik
			\tau_m v_{m0} -1 \pm \sqrt{G} }{2 \tau_m},~
		\dfrac{\dfrac{\partial p_c}{\partial \rho_c}ik \rho_{c0} \tau_c -2ik
			\tau_c v_{c0} -1 \pm \sqrt{G} }{2 \tau_c},
	\end{equation*} where
	\begin{equation*}
		G ={\left(\dfrac{\partial p_m}{\partial \rho_m}\right)^2}i^2k^2 \rho_{m0}^2 \tau_m^2 -4 \psi ik \rho_{m0} \tau_m v_{em}' - 2 \dfrac{\partial p_m}{\partial \rho_m}ik \rho_{m0} \tau_m +1.
	\end{equation*}
	For stability of the proposed model, changes in densities and velocities need to decrease as time progresses. Hence, the real parts $Re(r_{1 \pm, ~ 2 \pm}),$ of $det(J)$ must be strictly negative. That is to say; $Re(r_{1 \pm}), Re(r_{2 \pm}) < 0.$ However, the square root contains a complex number. So, it becomes difficult to determine the sign. Hence, we apply the formula \cite{helbing2013controversy}
	\begin{equation}
		\sqrt{R \pm i |I|} = \sqrt{\dfrac{1}{2}\left(\sqrt{R^2 + I^2} + R\right)} \pm i \sqrt{\dfrac{1}{2}\left(\sqrt{R^2+I^2} - R\right)},
		\label{RealPartOfw}
	\end{equation}
	where
	\begin{eqnarray*}
		R_c &=& 1 - {\left(\dfrac{\partial p_c(\rho_{c0})}{\partial \rho_c}\right)^2}k^2 \rho_{c0}^2 \tau_c^2,~
		I_c = 4 \psi_c k \rho_{c0} \tau_c v_{ec}' + 2 \dfrac{\partial p_c}{\partial \rho_c}k \rho_{c0} \tau_c, \\
		R_m &=& 1 - {\left(\dfrac{\partial p_m(\rho_{m0})}{\partial \rho_m}\right)^2}k^2 \rho_{m0}^2 \tau_m^2,~
		I_m = 4 \psi_m k \rho_{m0} \tau_m v_{em}' + 2 \dfrac{\partial p_m}{\partial \rho_m}k \rho_{m0} \tau_m. 
	\end{eqnarray*}
	Let us consider $Re(r_{1+}).$ To have $Re(r_{1+}) < 0,$ would mean
	\begin{equation*}
		-1 + \sqrt{\dfrac{1}{2} \left(\sqrt{R_m^2 + I_m^2} + R_m\right)} < 0.
	\end{equation*}
	After simplifying and substituting for $R_m$ and $I_m$ of (\ref{RealPartOfw}) into the previous the stability condition is derived as
	\begin{equation}
		\psi_m v_{em}' < - \dfrac{\partial p_m}{\partial \rho_{m}}.
		\label{stabilitycondition1}\end{equation}
	In the same way, the stability condition can be derived for $Re(r_{2+}) < 0$ and obtained as
	\begin{equation}
		\psi_c v_{em}' < - \dfrac{\partial p_c}{\partial \rho_{c}}.
		\label{stabilitycondition2}\end{equation}
	Data for the proposed model is chosen such that it satisfies the stability conditions in (\ref{stabilitycondition1}) and (\ref{stabilitycondition2}).
	
	\section{Numerical results}
	The following subsections describe vehicles' density, velocity and flow distributions when motorcycles' proportion is set to $20\%$ and $90\%.$ Considering (\ref{Rhom})and (\ref{Rhoc}), we look at freeway and congested traffic on roundabout. Roe's numerical scheme is utilized in which periodic boundary conditions are applied in order to represent a circular road. In an attempt to dissolve discontinuities at boundaries, an entropy fix is applied to the Roe's numerical scheme. Stability of the scheme is established by applying CFL stability condition. Table \ref{ParameterValues} shows parameter values that are used for simulations. Note that when $\delta$ takes on values $0$ and $1,$ results become unstable (undefined). Figures \ref{DensityVelocityFreeRoundaboutDelta20} - \ref{DensityVelocityCongestRoundabout3D} display simulation results of the proposed model.
	
	\begin{table}[!htp]
		\centering
		\caption{Simulation parameters.}
		\begin{tabular}{l|c|r}
			\hline
			\textbf{Description} & \textbf{Value} & Source\\ \hline
			Length of the road & 200m & \\ 
			Road step, h & 5 m & \cite{khan2019macroscopic}\\
			Time step & 0.05 s & \\
			Relaxation time,$\tau_m,~ \tau_c$ & 2s, 2.5s & \cite{khan2019macroscopic}\\
			Equilibrium velocity $v_e(\rho_m), v_e(\rho_c)$ & Greenshields & \cite{khan2019macroscopic}\\
			Maximum normalized density & $\rho=\rho_m+\rho_c = 1$ & \cite{khan2019macroscopic}\\
			Maximum speed, $v_{m\max}$ & 11m/s & \cite{mohan2017heterogeneous} \\
			Maximum speed, $v_{c\max}$ & 13.8m/s & \cite{Ministryofworksandtransport} \\
			Maximum area occupancy, $AO_m^{\max}$ & 0.85 & \cite{mohan2021multi} \\
			Maximum area occupancy, $AO_c^{\max}$ & 0.74 & \cite{mohan2021multi} \\		
			$\gamma_i,~ i=m,c$ & 2.23, 2.12 & \cite{mohan2021multi}\\
			Width of the road & 12m & \cite{Ministryofworksandtransport} \\
			Width of a car & 1.6m & \cite{arasan2008measuring}\\
			Vehicle class proportion, $\delta$ & 20\%, 90\% & \\
			Length of a car & 4m & \cite{arasan2008measuring}\\
			Length of a motorcycle & 1.8m & \cite{arasan2008measuring}\\
			Simulation time, T & 60 seconds & \cite{khan2019macroscopic} \\ \hline
		\end{tabular}
		\label{ParameterValues}
	\end{table}
	\subsection{Freeway traffic on a roundabout}
\begin{figure}[!htp]
		\vspace{0pt}
		\centering
		\begin{subfigure}[t]{0.3\textwidth}
			\vspace{0pt}
			\includegraphics[width=\linewidth]{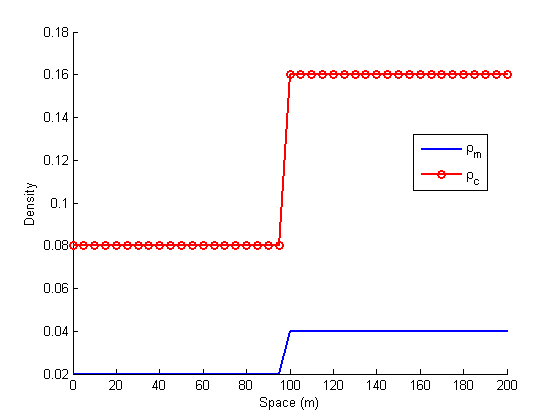}
			\caption{Density vs x at $0 s$}\label{FreewayDensityRoundaboutAt0s}
		\end{subfigure}
		\begin{subfigure}[t]{0.3\textwidth}
			\vspace{0pt}
			\includegraphics[width=\linewidth]{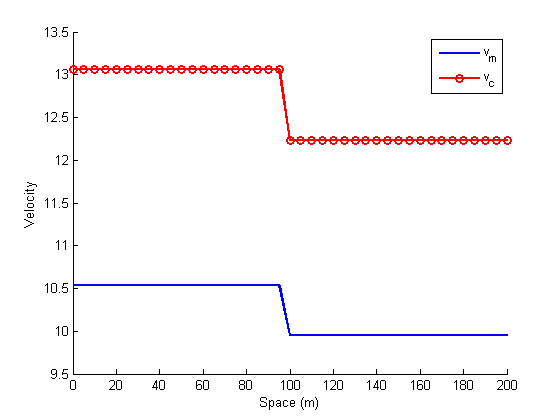}
			\caption{Velocity vs x at $0 s$}\label{FreewayVelocityRoundaboutAt0s}
		\end{subfigure}
		\begin{subfigure}[t]{0.3\textwidth}
			\vspace{0pt}
			\includegraphics[width=\linewidth]{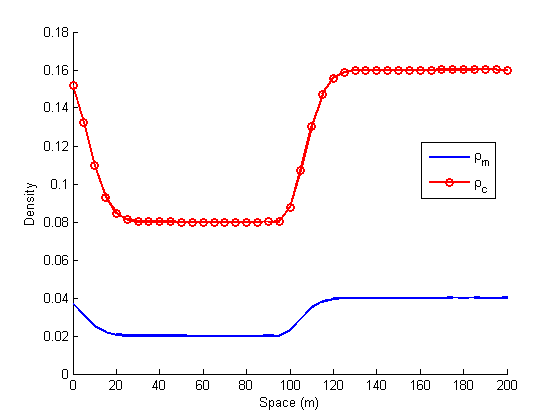}
			\caption{Density vs x at $1 s$}\label{FreewayDensityRoundaboutAt1s}
		\end{subfigure}
		\begin{subfigure}[t]{0.3\textwidth}
			\vspace{0pt}
			\includegraphics[width=\linewidth]{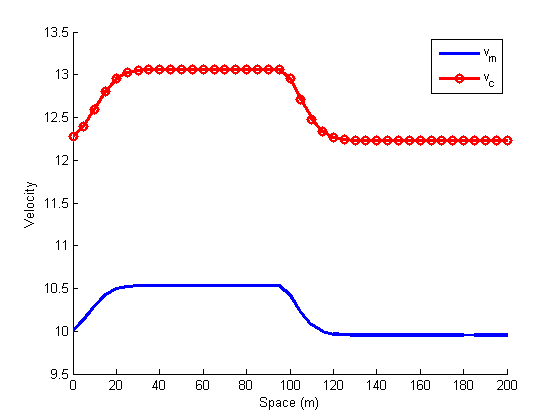}
			\caption{Velocity vs x at $1 s$ }\label{FreewayVelocityRoundaboutAt1s}
		\end{subfigure}
		\begin{subfigure}[t]{0.3\textwidth}
			\vspace{0pt}
			\includegraphics[width=\linewidth]{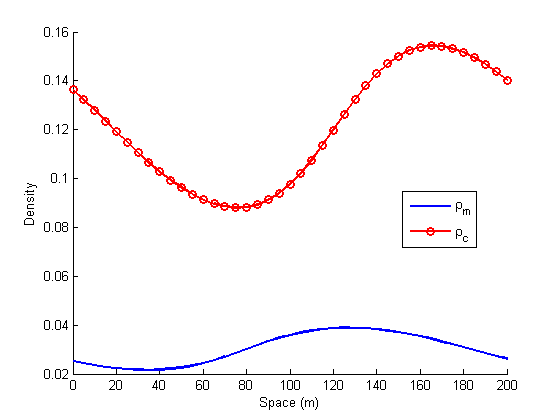}
			\caption{Density vs x at $20 s$}\label{FreewayDensityRoundaboutAt20s}
		\end{subfigure}
		\begin{subfigure}[t]{0.3\textwidth}
			\vspace{0pt}
			\includegraphics[width=\linewidth]{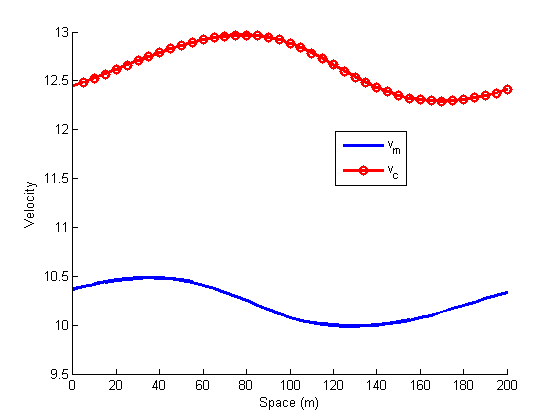}
			\caption{Velocity vs x at $20 s$}\label{FreewayVelocityRoundaboutAt20s}
		\end{subfigure}
		\begin{subfigure}[t]{0.3\textwidth}
			\vspace{0pt}
			\includegraphics[width=\linewidth]{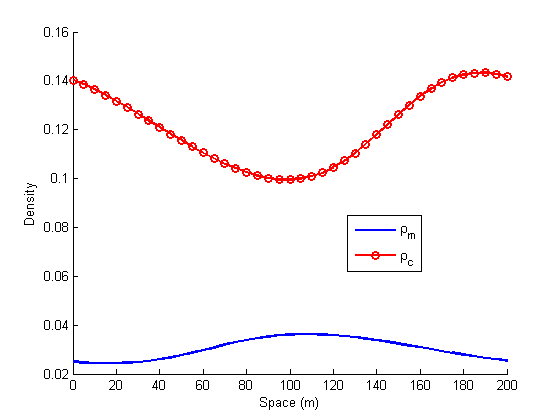}
			\caption{Density vs x at $40 s$}\label{FreewayDensityRoundaboutAt40s}
		\end{subfigure}
		\begin{subfigure}[t]{0.3\textwidth}
			\vspace{0pt}
			\includegraphics[width=\linewidth]{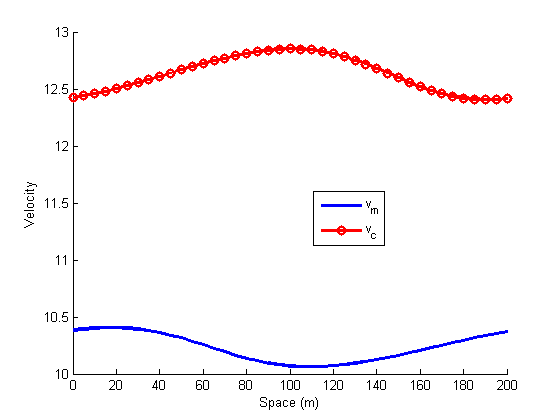}
			\caption{Velocity vs x, at $40 s$}\label{FreewayVelocityRoundaboutAt40s}
		\end{subfigure}
		\begin{subfigure}[t]{0.3\textwidth}
			\vspace{0pt}
			\includegraphics[width=\linewidth]{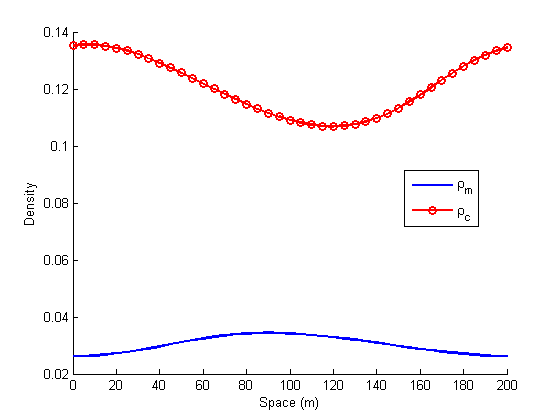}
			\caption{Density vs x at $60 s$}\label{FreewayDensityRoundaboutAt60s}
		\end{subfigure}
		\begin{subfigure}[t]{0.3\textwidth}
			\vspace{0pt}
			\includegraphics[width=\linewidth]{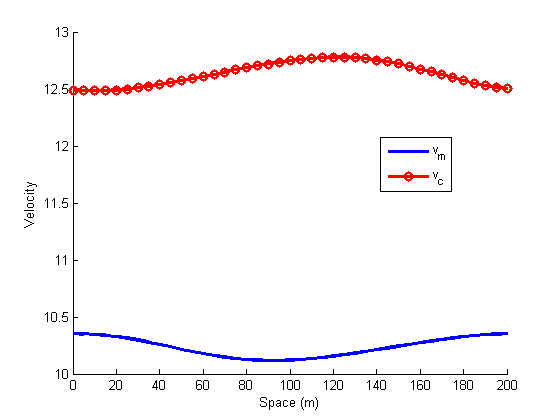}
			\caption{Velocity vs x, at $60 s$}\label{FreewayVelocityRoundaboutAt60s}
		\end{subfigure}
		
\caption{Densities and velocities of the proposed model on a $200 m$ freeway circular road when $\delta = 20\%,$ at $T = 0s, 1 s, 20s, 40s$ and $60 s.$  }	\label{DensityVelocityFreeRoundaboutDelta20}
		\vspace{-0.1cm}
	\end{figure}

\begin{figure}
		\vspace{0pt}
		\centering 
		\begin{subfigure}[t]{0.32\textwidth}
			\vspace{0pt}
			\includegraphics[width=\linewidth]{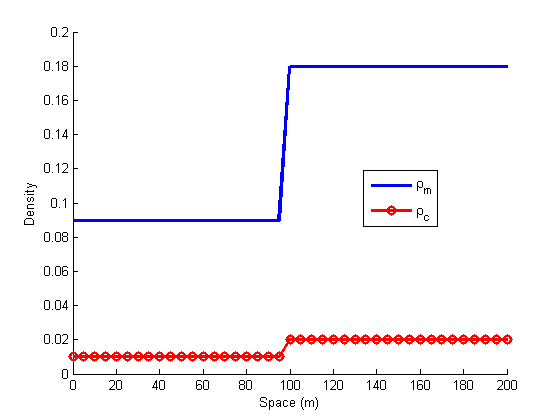}
			\caption{Density vs x at $0 s$} \label{FreewayDensityRoundaboutAt0sDelta90}
		\end{subfigure}
		\begin{subfigure}[t]{0.32\textwidth}
			\vspace{0pt}
			\includegraphics[width=\linewidth]{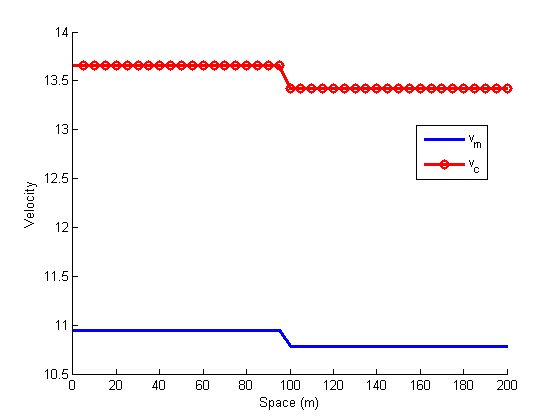}
			\caption{Velocity vs x at $0 s$}\label{FreewayVelocityRoundaboutAt0sDelta90}
		\end{subfigure}
		\begin{subfigure}[t]{0.32\textwidth}
			\vspace{0pt}
			\includegraphics[width=\linewidth]{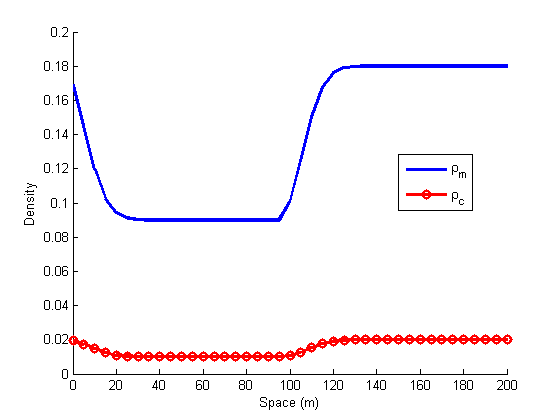}
			\caption{Density vs x at $1 s$}\label{FreewayDensityRoundaboutAt1sDelta90}
		\end{subfigure}
		\begin{subfigure}[t]{0.32\textwidth}
			\vspace{0pt}
			\includegraphics[width=\linewidth]{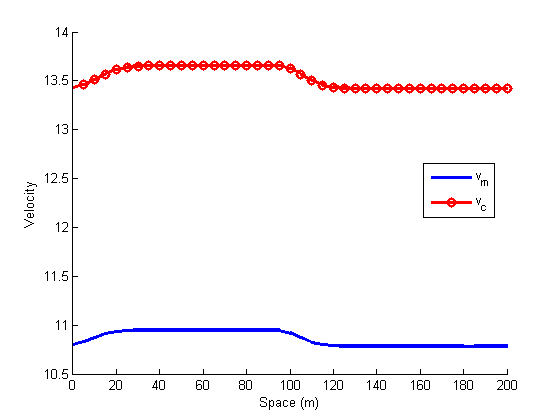}
			\caption{Velocity vs x at $1 s$}\label{FreewayVelocityRoundaboutAt1sDelta90}
		\end{subfigure}
		\begin{subfigure}[t]{0.32\textwidth}
			\vspace{0pt}
			\includegraphics[width=\linewidth]{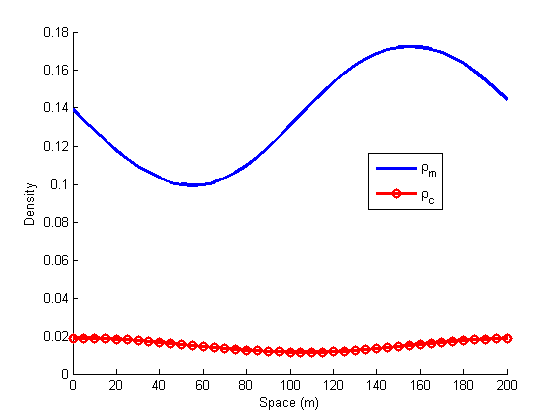}
			\caption{Density vs x, at $20 s$}\label{FreewayDensityRoundaboutAt20sDelta90}
		\end{subfigure}
		\begin{subfigure}[t]{0.32\textwidth}
			\vspace{0pt}
			\includegraphics[width=\linewidth]{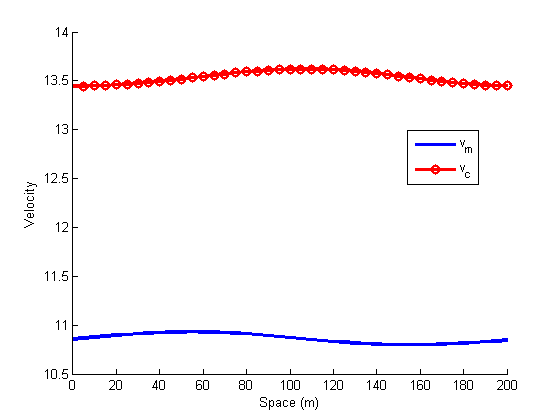}
			\caption{Velocity vs x at $20 s$}\label{FreewayVelocityRoundaboutAt20sDelta90}
		\end{subfigure}
		\begin{subfigure}[t]{0.32\textwidth}
			\vspace{0pt}
			\includegraphics[width=\linewidth]{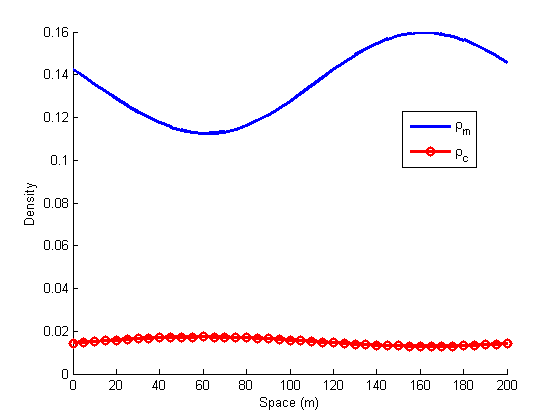}
			\caption{Density vs x at $40 s$}\label{FreewayDensityRoundaboutAt40sDelta90}
		\end{subfigure}
		\begin{subfigure}[t]{0.32\textwidth}
			\vspace{0pt}
			\includegraphics[width=\linewidth]{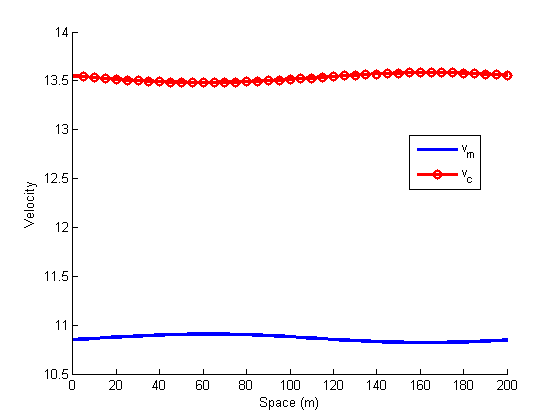}
			\caption{Velocity vs x at $40 s$}\label{FreewayVelocityRoundaboutAt40sDelta90}
		\end{subfigure}
		\begin{subfigure}[t]{0.32\textwidth}
			\vspace{0pt}
			\includegraphics[width=\linewidth]{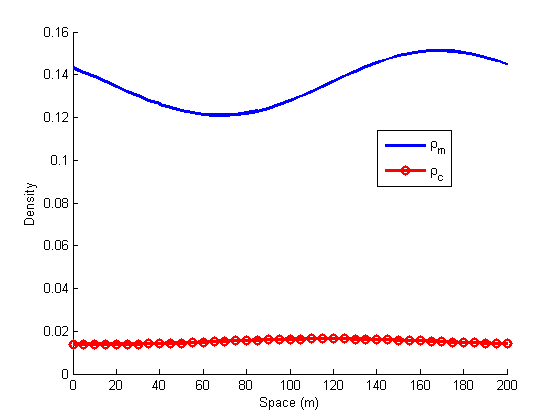}
			\caption{Density vs x at $60 s$}\label{FreewayDensityRoundaboutAt60sDelta90}
		\end{subfigure}
		\begin{subfigure}[t]{0.32\textwidth}
			\vspace{0pt}
			\includegraphics[width=\linewidth]{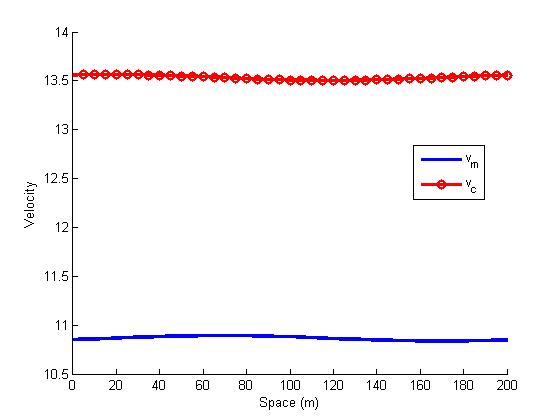}
			\caption{Velocity vs x at $60 s$}\label{FreewayVelocityRoundaboutAt60sDelta90}
		\end{subfigure}
		\caption{Densities and velocities of the proposed model on a $200 m$ freeway circular road when $\delta = 90\%,$ at $T = 0s, 1s, 20s, 40 s$ and $60 s.$ }	\label{DensityVelocityFreeRoundaboutDelta90}
		
		\vspace{-0.1 cm}
		
	\end{figure}

	Under this subsection, we consider a free traffic flow with motorcycles proportion $\delta,$ initially fixed at $20\%,$ for $0s, 1 s, ~ 20 s, ~ 40 s$ and $60 s.$ The initial total density of traffic flow denoted $\rho_0,$ is set to be 
	\begin{eqnarray}
		\rho_0 &=&
		\begin{cases}
			0.1, &\text{for}~ x < 100 \\
			0.2,& \text{for}~ x \geq 100.
		\end{cases} \label{FreeRoundabout}
	\end{eqnarray}
	Since $\delta = 20\%$ it follows that $\rho_m = 0.2 \rho_0$ and $\rho_c = 0.8 \rho_0.$ The initial density distributions of the respective vehicle classes are plotted in Figure \ref{FreewayDensityRoundaboutAt0s}, and corresponding velocities depicted by Figure \ref{FreewayVelocityRoundaboutAt0s}. In Figure \ref{DensityVelocityFreeRoundaboutDelta20}, a shock wave develops at $1 s$ for both vehicle classes at $100 m$ (when the low density stream meets high density stream). After $20s,$ the shock wave is smoothed and both vehicle classes move with velocities close to their respective maximum values. That is; $13.8 m/s$ and $11 m/s$ for cars and motorcycles, respectively (see Figures \ref{FreewayVelocityRoundaboutAt20s} - \ref{FreewayVelocityRoundaboutAt60s}). Next, $\delta$ is set to $90\%,$ such that $\rho_m = 0.9\rho_0$ and $\rho_c = 0.1\rho_0,$ for motorcycles and cars, respectively. The plot of such initial density profile is shown in Figure \ref{FreewayDensityRoundaboutAt0sDelta90} with corresponding velocities in Figure \ref{FreewayVelocityRoundaboutAt0sDelta90}. Same results are observed when $\delta=90\%$ (see Figures \ref{FreewayVelocityRoundaboutAt20sDelta90} - \ref{FreewayVelocityRoundaboutAt60sDelta90}) Therefore, it can be concluded that there is minimal effect of motorcycles proportion during freeway traffic flow.
	
	\subsection{Congested traffic on a roundabout}
\begin{figure}[!htp]
		
		\vspace{0pt}
		
		\centering
		\begin{subfigure}[t]{0.32\textwidth}
			\vspace{0pt}
			\includegraphics[width=\linewidth]{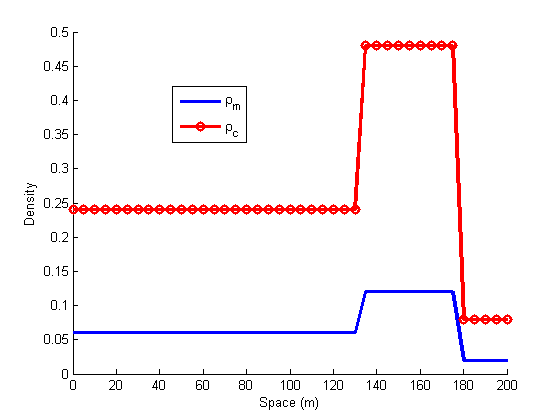}
			\caption{Density vs x at $0 s$}\label{CongestDensityRoundaboutAt0sDelta20}
		\end{subfigure}
		\begin{subfigure}[t]{0.32\textwidth}
			\vspace{0pt}
			\includegraphics[width=\linewidth]{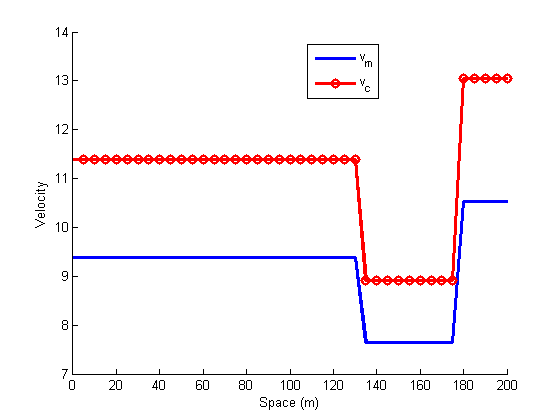}
			\caption{Velocity vs x at $0 s$}\label{CongestVelocityRoundaboutAt0sDelta20}
		\end{subfigure}
		\begin{subfigure}[t]{0.32\textwidth}
			\vspace{0pt}
			\includegraphics[width=\linewidth]{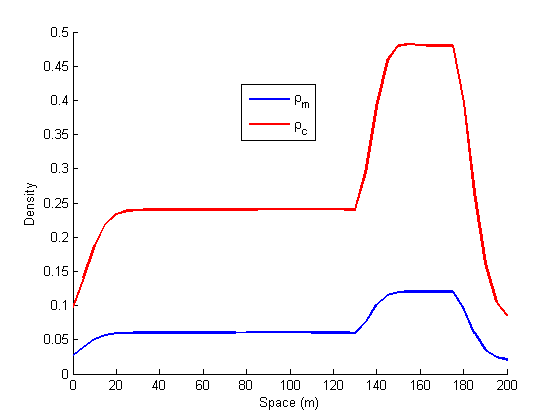}
			\caption{Density vs x at $1 s$}\label{CongestDensityRoundaboutAt1s}
		\end{subfigure}
		\begin{subfigure}[t]{0.32\textwidth}
			\vspace{0pt}
			\includegraphics[width=\linewidth]{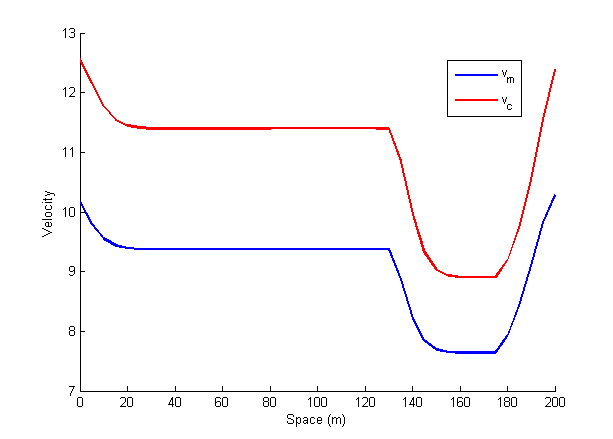}
			\caption{Velocity vs x at $1 s$}\label{CongestVelocityRoundaboutAt1s}
		\end{subfigure}
		\begin{subfigure}[t]{0.32\textwidth}
			\vspace{0pt}
			\includegraphics[width=\linewidth]{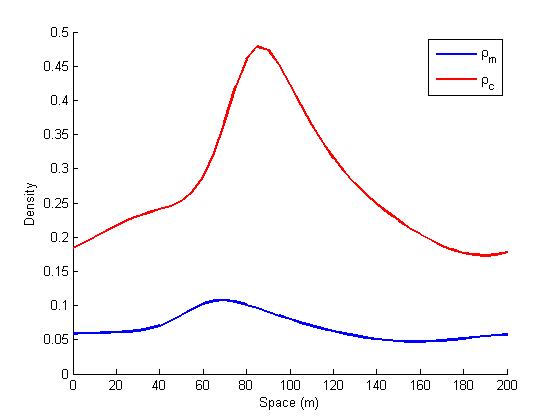}
			\caption{Density vs x at $20 s$}\label{CongestDensityRoundaboutAt20s}
		\end{subfigure}
		\begin{subfigure}[t]{0.32\textwidth}
			\vspace{0pt}
			\includegraphics[width=\linewidth]{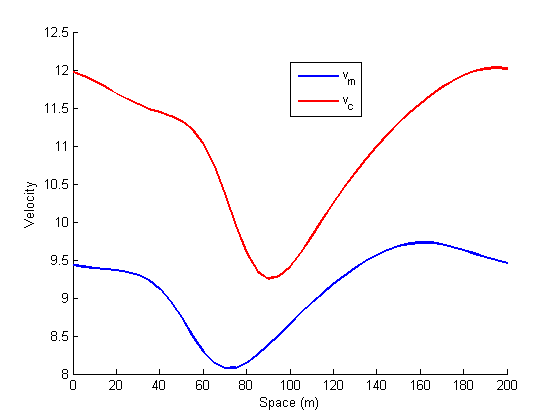}
			\caption{Velocity vs x at $20 s$ }\label{CongestVelocityRoundaboutAt20s}
		\end{subfigure}
		\begin{subfigure}[t]{0.32\textwidth}
			\vspace{0pt}
			\includegraphics[width=\linewidth]{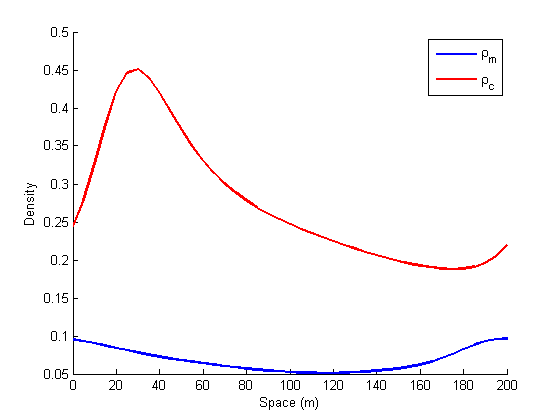}
			\caption{Density vs x at $40 s$}\label{CongestDensityRoundaboutAt40s}
		\end{subfigure}
		\begin{subfigure}[t]{0.32\linewidth}
			\vspace{0pt}
			\includegraphics[width=\linewidth]{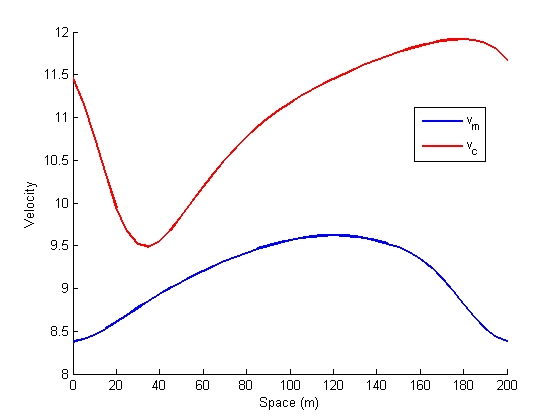}
			\caption{Velocity vs x, at $40 s$}\label{CongestVelocityRoundaboutAt40s}
		\end{subfigure}
		\begin{subfigure}[t]{0.32\textwidth}
			\vspace{0pt}
			\includegraphics[width=\linewidth]{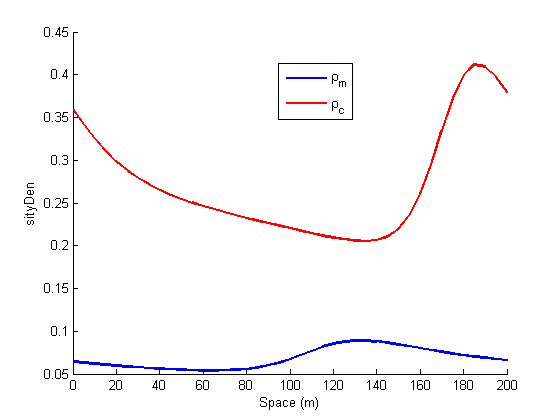}
			\caption{Density vs x at $60 s$}\label{CongestDensityRoundaboutAt60s}
		\end{subfigure}
		\begin{subfigure}[t]{0.32\textwidth}
			\vspace{0pt}
			\includegraphics[width=\linewidth]{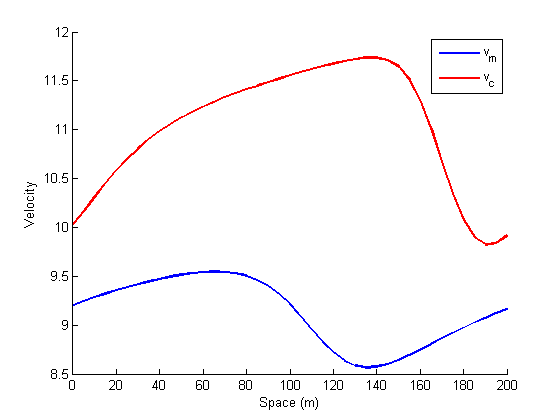}
			\caption{Velocity vs x at $60 s$ }\label{CongestVelocityRoundaboutAt60s}
		\end{subfigure}
		\caption{Densities and velocities of the proposed model on a $200 m$ congested circular road when $\delta = 20\%,$ at time $T = 0s, 1s, 20s, 40 s$ and $60 s.$}	\label{DensityVelocityCongestRoundaboutDelta20}
		
		\vspace{-0.1cm}
		
	\end{figure}

	Now we consider a congested $200 m$ circular road with total initial density, $\rho_0$ given by 
	\begin{eqnarray}
		\rho_0 &=&
		\begin{cases}
			0.3, &\text{if}~ x \leq 130 \\
			0.6,& \text{if}~ 130 < x < 180 \\
			0.1, & \text{if}~ x \geq 180.
		\end{cases}
	\end{eqnarray}
	The proposed model densities and velocities at time, $T=0s, 1s,~ 20 s,~ 40s$ and $60s,$ with $\delta$ initially fixed at $20\%$ are displayed in Figure \ref{DensityVelocityCongestRoundaboutDelta20}. Plots of initial densities and velocities are displayed by Figures \ref{CongestDensityRoundaboutAt0sDelta20} and \ref{CongestVelocityRoundaboutAt0sDelta20}, respectively. It is observed in Figures \ref{CongestDensityRoundaboutAt1s} - \ref{CongestVelocityRoundaboutAt60s} that at all selected times, both vehicle classes develop shock waves as high and low densities of vehicles meet. However, the shock waves are smoothed at $20s$ and beyond, since an entropy fix is applied to Roe numerical scheme. We further observe that velocities for both vehicle classes decrease far away from their maximum values as time progresses. Different from the freeway traffic flow scenario, velocities of both vehicle classes are non uniform at all times since they are congested. Cars form a rarefaction wave at $T=1s$ and $180m,$  but propagates to the left at $20s$ and $40s$ and later disappears at $60s.$ A rarefaction wave develops for motorcycles at $T=1s$ and $180m$ but disappears at $20s.$ A shock wave develops for both vehicle classes at $T=1s$ and $130m,$ but for cars it propagates to the left at $20s$ and $40s$ and later moves back to the right at $60s.$ At $20s$ and $60s,$ motorcycles hit jam before cars. At $40s,$ cars hit jam faster than motorcycles. At any time, jam is created in one region for cars while motorcycles density is maintained at low values.
	
\begin{figure}[!htp]
		\vspace{0pt}
		\centering
		\begin{subfigure}[t]{0.32\textwidth}
			\vspace{0pt}
			\includegraphics[width=\linewidth]{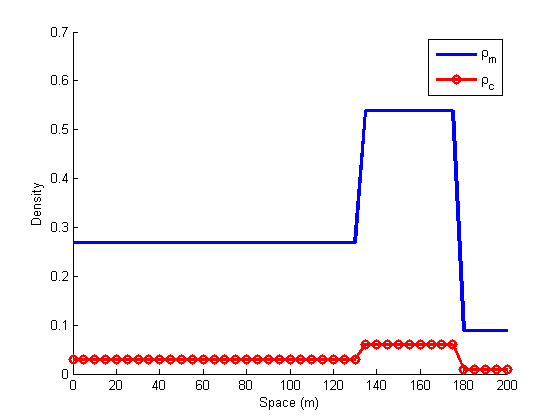}
			\caption{Density vs x at $0 s$}\label{CongestDensityRoundaboutAt0sDelta90}
		\end{subfigure}
		\begin{subfigure}[t]{0.32\textwidth}
			\vspace{0pt}
			\includegraphics[width=\linewidth]{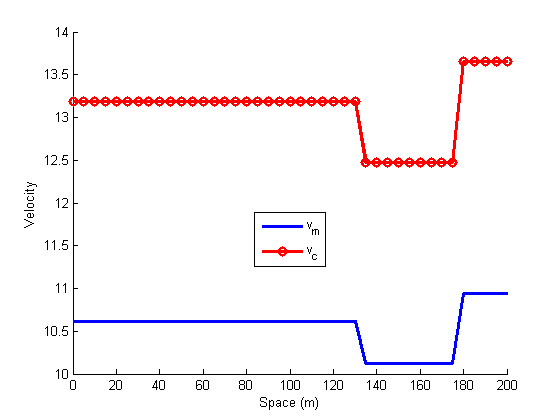}
			\caption{Velocity vs x at $0 s$}\label{CongestVelocityRoundaboutAt0sDelta90}
		\end{subfigure}
		\begin{subfigure}[t]{0.32\textwidth}
			\vspace{0pt}
			\includegraphics[width=\linewidth]{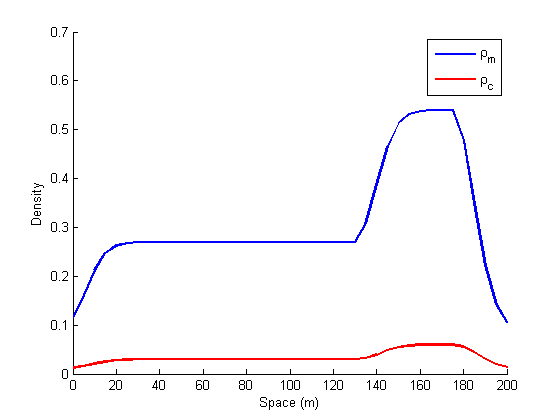}
			\caption{Density vs x at $1 s$}\label{CongestDensityRoundaboutAt1sDelta90}
		\end{subfigure}
		\begin{subfigure}[t]{0.32\textwidth}
			\vspace{0pt}
			\includegraphics[width=\linewidth]{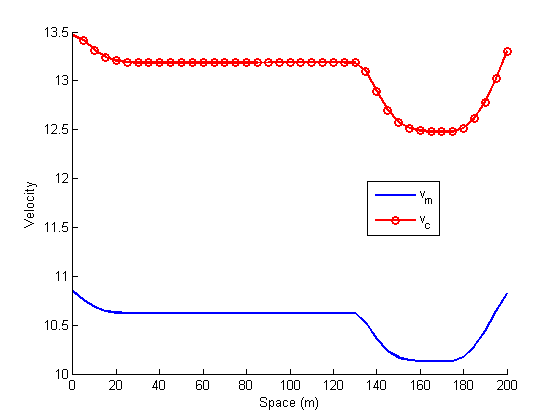}
			\caption{Velocity vs x at $1 s$}\label{CongestVelocityRoundaboutAt1sDelta90}
		\end{subfigure}
		\begin{subfigure}[t]{0.32\textwidth}
			\vspace{0pt}
			\includegraphics[width=\linewidth]{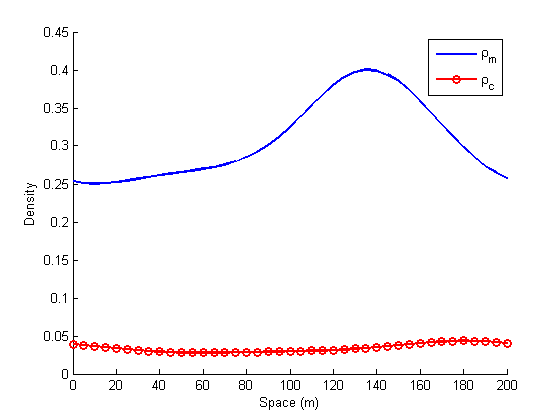}
			\caption{Density vs x at $20 s$}\label{CongestDensityRoundaboutAt20sDelta90}
		\end{subfigure}
		\begin{subfigure}[t]{0.32\textwidth}
			\vspace{0pt}
			\includegraphics[width=\linewidth]{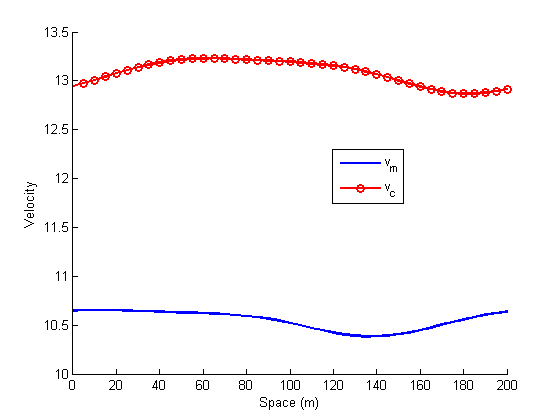}
			\caption{Velocity vs x at $20 s$ }\label{CongestVelocityRoundaboutAt20sDelta90}
		\end{subfigure}
		\begin{subfigure}[t]{0.32\textwidth}
			\vspace{0pt}
			\includegraphics[width=\linewidth]{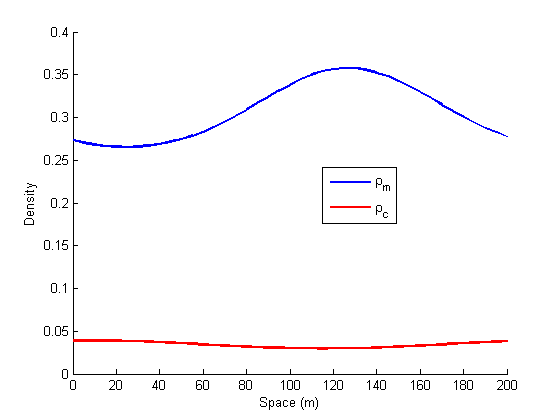}
			\caption{Density vs x at $40 s$}\label{CongestDensityRoundaboutAt40sDelta90}
		\end{subfigure}
		\begin{subfigure}[t]{0.32\linewidth}
			\vspace{0pt}
			\includegraphics[width=\linewidth]{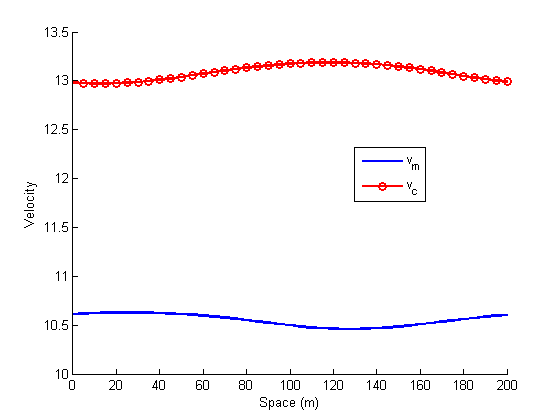}
			\caption{Velocity vs x, at $40 s$}\label{CongestVelocityRoundaboutAt40sDelta90}
		\end{subfigure}
		\begin{subfigure}[t]{0.32\textwidth}
			\vspace{0pt}
			\includegraphics[width=\linewidth]{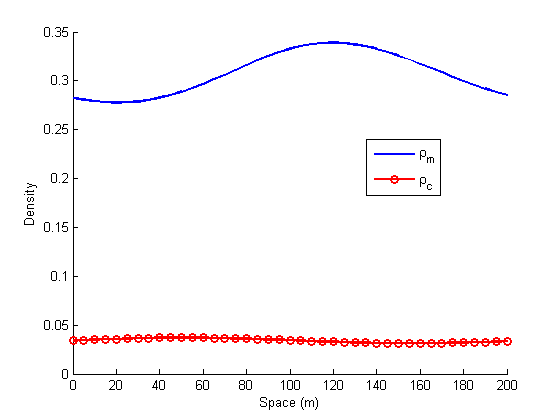}
			\caption{Density vs x at $60 s$}\label{CongestDensityRoundaboutAt60sDelta90}
		\end{subfigure}
		\begin{subfigure}[t]{0.32\textwidth}
			\vspace{0pt}
			\includegraphics[width=\linewidth]{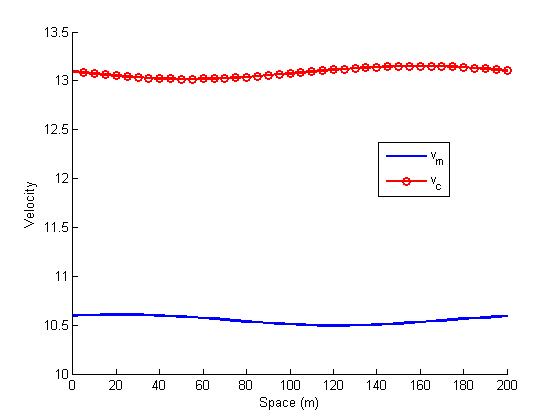}
			\caption{Velocity vs x at $60 s$ }\label{CongestVelocityRoundaboutAt60sDelta90}
		\end{subfigure}
		\caption{Densities and velocities of the proposed model on a $200 m$ congested circular road when $\delta = 90\%,$ at $0s, 1s, 20s, 40 s$ and $60 s.$}	\label{CongestDensityVelocityRoundaboutDelta90}
		
		\vspace{-0.1cm}
		
	\end{figure}

\begin{figure}[H]
		\vspace{0pt}
		\centering
		\begin{subfigure}[t]{0.32\textwidth}
			\vspace{0pt}
			\includegraphics[width=\linewidth]{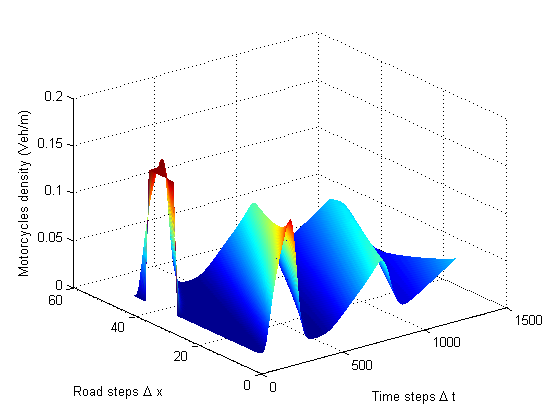}
			\caption{$\rho_m$ vs x vs t, $\delta=20\%$}\label{CongestDensitymRoundaboutAtDelta203D}
		\end{subfigure}
		\begin{subfigure}[t]{0.32\textwidth}
			\vspace{0pt}
			\includegraphics[width=\linewidth]{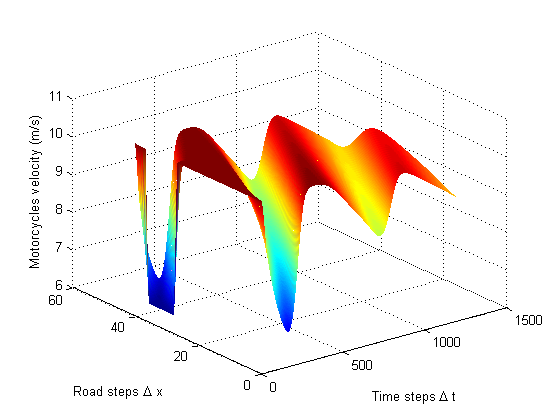}
			\caption{$v_m$ vs x vs t, $\delta=20\%$}
			\label{CongestVelocitymRoundaboutAtDelta203D}
		\end{subfigure}
		\begin{subfigure}[t]{0.32\textwidth}
			\vspace{0pt}
			\includegraphics[width=\linewidth]{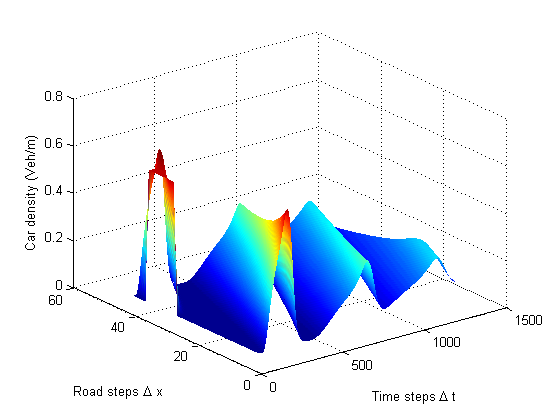}
			\caption{$\rho_c$ vs x vs t, at $\delta=20\%$}\label{CongestDensitycRoundaboutAtDelta203D}
		\end{subfigure}
		\begin{subfigure}[t]{0.32\textwidth}
			\vspace{0pt}
			\includegraphics[width=\linewidth]{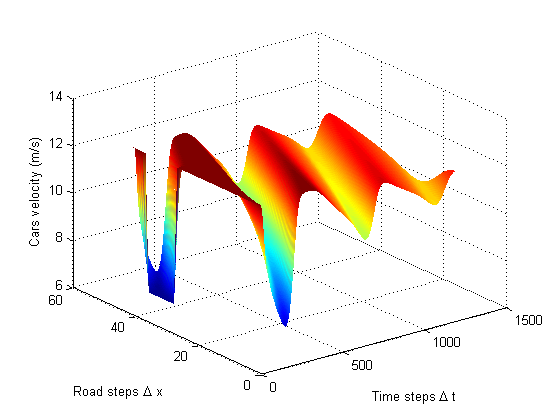}
			\caption{$v_c$ vs x vs t, $\delta=20\%$}
			\label{CongestVelocitycRoundaboutAtDelta203D}
		\end{subfigure}
		\begin{subfigure}[t]{0.32\textwidth}
			\vspace{0pt}
			\includegraphics[width=\linewidth]{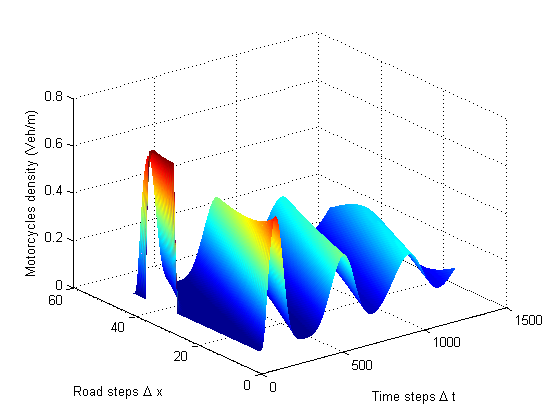}
			\caption{$\rho_m$ vs x vs t, $\delta =90\%$}\label{CongestDensitymRoundaboutAtDelat903D}
		\end{subfigure}
		\begin{subfigure}[t]{0.32\textwidth}
			\vspace{0pt}
			\includegraphics[width=\linewidth]{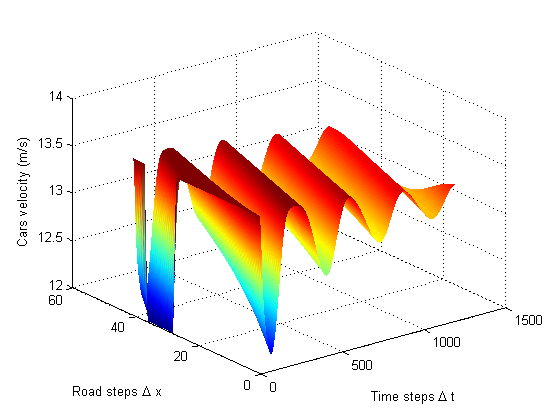}
			\caption{$v_m$ vs x vs t, $\delta =90\%$}\label{CongestVelocitymRoundaboutAtDelat903D}
		\end{subfigure}
		\begin{subfigure}[t]{0.32\textwidth}
			\vspace{0pt}
			\includegraphics[width=\linewidth]{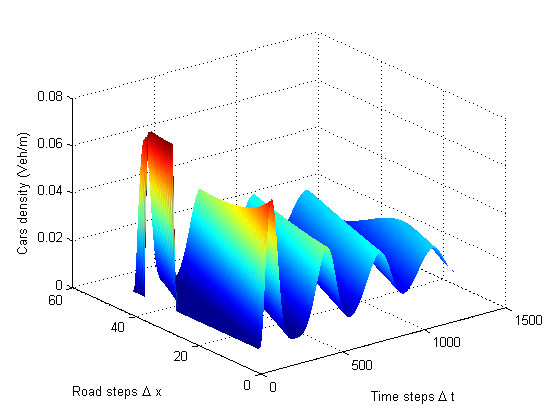}
			\caption{$\rho_c$ vs x vs t, $\delta =90\%$}\label{CongestDensitycRoundaboutAtDelat903D}
		\end{subfigure}
		\begin{subfigure}[t]{0.32\textwidth}
			\vspace{0pt}
			\includegraphics[width=\linewidth]{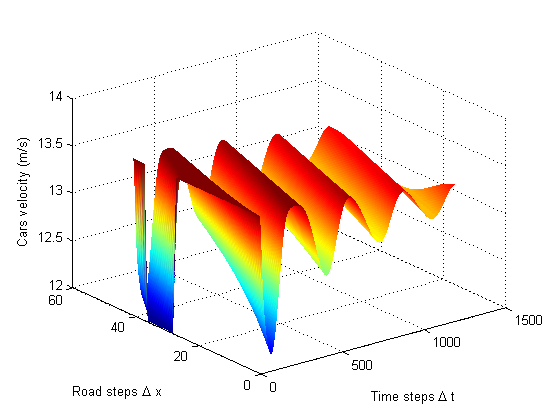}
			\caption{$v_c$ vs x vs t, $\delta=90\%$}
			\label{CongestVelocitycRoundaboutAtDelat903D}
		\end{subfigure}
		\caption{Densities and velocities of the proposed model on a $200 m$ congested circular road when $\delta = 20\%$ and $90\%,$ at $60 s.$ }	\label{DensityVelocityCongestRoundabout3D}
		
		\vspace{-0.1cm}
		
	\end{figure}

	Next $\delta$ is set to be $90\%,$ that is; fewer cars and many motorcycles. The initial density and velocity distributions are indicated in Figures \ref{CongestDensityRoundaboutAt0sDelta90} and \ref{CongestVelocityRoundaboutAt0sDelta90}. Density for cars is observed to be almost uniform after $1s$ while for motorcycles is non uniform at all times. A shock wave is shown to clear faster for cars than for motorcycles. A stationary shock wave but with decreasing amplitude develops for motorcycles in one region after $1s,$ where they find difficulties squeezing through spaces unfilled by the many cars. Figures \ref{CongestDensitymRoundaboutAtDelta203D} - \ref{CongestVelocitycRoundaboutAtDelat903D} display the space-time evolution of densities and velocities for both vehicle classes, when $\delta = 20\%$ and $90\%.$ The densities and velocities become smooth as time progresses. 
	
	\section{Conclusion} 
	The aim of this paper is to determine the effect of proportional densities on the flow of traffic. A new expression for area occupancy, defined in terms of proportion densities is proposed. Roe's scheme is explicitly applied, where Roe matrix, averaged velocity, averaged pressure, wave strengths are derived. An entropy condition is applied to Roe decomposition. Stability conditions for each vehicle class are obtained by linear stability analysis whereas stability of the numerical scheme is determined by the CFL condition.
	
	Qualitative analysis shows that Roe decomposition scheme can be applied to AR model and other variants. Numerical results show that when motorcycles are more than cars, all vehicle classes move with higher velocities than when motorcycles are fewer. A shock wave is observed to clear faster for cars than for motorcycles since motorcycles are slower than cars. On the other hand if motorcycles are fewer than cars, both vehicle classes develop shock waves but this effect is more felt during congested traffic flow.
	Traffic jam is shown to occur much earlier if cars dominate the road than when motorcycles do. At any time, whenever cars develop shock wave motorcycles too develop the shock. Densities and velocities for each vehicle class become smoother over time and were shown to remain within limits. Thus, the results obtained from the proposed model are realistic. Total traffic flow increases with increase of motorcycles proportion. 
	As a control for traffic jam, enforcement of lane discipline is recommended.
	\section*{Acknowledgment}
	The first author acknowledges her other research supervisors Prof. J.Y.T Mugisha of Department of Mathematics, Makerere University, P.O Box 7062, Kampala, Uganda and Prof. Semu Mitiku Kassa of Department of Mathematics and Statistical Sciences, Botswana
	International University of Science and Technology (BIUST), P/Bag $16,$ Palapye, Botswana for providing technical guidance and mentorship. The first author also acknowledges the Sida bilateral program with Makerere University, 2015-2020, project 316 ``Capacity building in Mathematics and its applications'' for providing financial support.

\end{document}